%
%
%
%



\magnification 1200
\hsize 13.2cm
\vsize 20cm
\parskip 3pt plus 1pt
\parindent 5mm

\def\\{\hfil\break}


\font\seventeenbf=cmbx10 at 17.28pt

\font\twelvebf=cmbx10 at 12pt
\font\eightbf=cmbx8
\font\sixbf=cmbx6

\font\eighti=cmmi8
\font\sixi=cmmi6

\font\eightrm=cmr8
\font\sixrm=cmr6

\font\eightsy=cmsy8
\font\sixsy=cmsy6

\font\eightit=cmti8
\font\eighttt=cmtt8
\font\eightsl=cmsl8

\font\seventeenbsy=cmbsy10 at 17.28pt

\font\twelvebsy=cmbsy10 at 12pt
\font\tenbsy=cmbsy10
\font\eightbsy=cmbsy8
\font\sevenbsy=cmbsy7
\font\sixbsy=cmbsy6
\font\fivebsy=cmbsy5

\font\tenmsa=msam10

\font\sevenmsa=msam7
\font\fivemsa=msam5
\newfam\msafam
  \textfont\msafam=\tenmsa
  \scriptfont\msafam=\sevenmsa
  \scriptscriptfont\msafam=\fivemsa

\font\tenmsb=msbm10
\font\eightmsb=msbm8
\font\sevenmsb=msbm7
\font\fivemsb=msbm5
\newfam\msbfam
  \textfont\msbfam=\tenmsb
  \scriptfont\msbfam=\sevenmsb
  \scriptscriptfont\msbfam=\fivemsb
\def\Bbb{\fam\msbfam\tenmsb}

\font\tenCal=eusm10
\font\sevenCal=eusm7
\font\fiveCal=eusm5
\newfam\Calfam
  \textfont\Calfam=\tenCal
  \scriptfont\Calfam=\sevenCal
  \scriptscriptfont\Calfam=\fiveCal
\def\Cal{\fam\Calfam\tenCal}

\font\teneuf=eusm10
\font\teneuf=eufm10
\font\seveneuf=eufm7
\font\fiveeuf=eufm5
\newfam\euffam
  \textfont\euffam=\teneuf
  \scriptfont\euffam=\seveneuf
  \scriptscriptfont\euffam=\fiveeuf

\font\seventeenbfit=cmmib10 at 17.28pt

\font\twelvebfit=cmmib10 at 12pt
\font\tenbfit=cmmib10
\font\eightbfit=cmmib8
\font\sevenbfit=cmmib7
\font\sixbfit=cmmib6
\font\fivebfit=cmmib5
\newfam\bfitfam
  \textfont\bfitfam=\tenbfit
  \scriptfont\bfitfam=\sevenbfit
  \scriptscriptfont\bfitfam=\fivebfit


\catcode`\@=11
\def\eightpoint{%
  \textfont0=\eightrm \scriptfont0=\sixrm \scriptscriptfont0=\fiverm
  \def\rm{\fam\z@\eightrm}%
  \textfont1=\eighti \scriptfont1=\sixi \scriptscriptfont1=\fivei
  \def\oldstyle{\fam\@ne\eighti}%
  \textfont2=\eightsy \scriptfont2=\sixsy \scriptscriptfont2=\fivesy
  \textfont\itfam=\eightit
  \def\it{\fam\itfam\eightit}%
  \textfont\slfam=\eightsl
  \def\sl{\fam\slfam\eightsl}%
  \textfont\bffam=\eightbf \scriptfont\bffam=\sixbf
  \scriptscriptfont\bffam=\fivebf
  \def\bf{\fam\bffam\eightbf}%
  \textfont\ttfam=\eighttt
  \def\tt{\fam\ttfam\eighttt}%
  \textfont\msbfam=\eightmsb
  \def\Bbb{\fam\msbfam\eightmsb}%
  \abovedisplayskip=9pt plus 2pt minus 6pt
  \abovedisplayshortskip=0pt plus 2pt
  \belowdisplayskip=9pt plus 2pt minus 6pt
  \belowdisplayshortskip=5pt plus 2pt minus 3pt
  \smallskipamount=2pt plus 1pt minus 1pt
  \medskipamount=4pt plus 2pt minus 1pt
  \bigskipamount=9pt plus 3pt minus 3pt
  \normalbaselineskip=9pt
  \setbox\strutbox=\hbox{\vrule height7pt depth2pt width0pt}%
  \let\bigf@ntpc=\eightrm \let\smallf@ntpc=\sixrm
  \normalbaselines\rm}
\catcode`\@=12

\def\eightpointbf{%
 \textfont0=\eightbf   \scriptfont0=\sixbf   \scriptscriptfont0=\fivebf
 \textfont1=\eightbfit \scriptfont1=\sixbfit \scriptscriptfont1=\fivebfit
 \textfont2=\eightbsy  \scriptfont2=\sixbsy  \scriptscriptfont2=\fivebsy
 \eightbf
 \baselineskip=10pt}

\def\tenpointbf{%
 \textfont0=\tenbf   \scriptfont0=\sevenbf   \scriptscriptfont0=\fivebf
 \textfont1=\tenbfit \scriptfont1=\sevenbfit \scriptscriptfont1=\fivebfit
 \textfont2=\tenbsy  \scriptfont2=\sevenbsy  \scriptscriptfont2=\fivebsy
 \tenbf}

\def\twelvepointbf{%
 \textfont0=\twelvebf   \scriptfont0=\eightbf   \scriptscriptfont0=\sixbf
 \textfont1=\twelvebfit \scriptfont1=\eightbfit \scriptscriptfont1=\sixbfit
 \textfont2=\twelvebsy  \scriptfont2=\eightbsy  \scriptscriptfont2=\sixbsy
 \twelvebf
 \baselineskip=14.4pt}

\def\seventeenpointbf{%
 \textfont0=\seventeenbf  \scriptfont0=\twelvebf  \scriptscriptfont0=\eightbf
 \textfont1=\seventeenbfit\scriptfont1=\twelvebfit\scriptscriptfont1=\eightbfit
 \textfont2=\seventeenbsy \scriptfont2=\twelvebsy \scriptscriptfont2=\eightbsy
 \seventeenbf
 \baselineskip=20.736pt}


\newdimen\srdim \srdim=\hsize
\newdimen\irdim \irdim=\hsize
\def\NOSECTREF#1{\noindent\hbox to \srdim{\null\dotfill ???(#1)}}
\def\SECTREF#1{\noindent\hbox to \srdim{\csname REF\romannumeral#1\endcsname}}
\def\INDREF#1{\noindent\hbox to \irdim{\csname IND\romannumeral#1\endcsname}}
\newlinechar=`\^^J
\def\openauxfile{
  \immediate\openin1\jobname.aux
  \ifeof1
  \message{^^JCAUTION\string: you MUST run TeX a second time^^J}
  \let\sectref=\NOSECTREF \let\indref=\NOSECTREF
  \else
  \input \jobname.aux
  \message{^^JCAUTION\string: if the file has just been modified you may
    have to run TeX twice^^J}
  \let\sectref=\SECTREF \let\indref=\INDREF
  \fi
  \message{to get correct page numbers displayed in Contents or Index
    Tables^^J}
  \immediate\openout1=\jobname.aux
  \let\END=\end \def\end{\immediate\closeout1\END}}

\newbox\titlebox   \setbox\titlebox\hbox{\hfil}
\newbox\sectionbox \setbox\sectionbox\hbox{\hfil}
\def\folio{\ifnum\pageno=1 \hfil \else \ifodd\pageno
           \hfil {\eightpoint\copy\sectionbox\kern8mm\number\pageno}\else
           {\eightpoint\number\pageno\kern8mm\copy\titlebox}\hfil \fi\fi}
\footline={\hfil}
\headline={\folio}

\def\titlerunning#1{\setbox\titlebox\hbox{\eightpoint #1}}
\def\title#1{\noindent\hfil$\smash{\hbox{\seventeenpointbf #1}}$\hfil
             \titlerunning{#1}\medskip}

\newcount\numbersection \numbersection=-1
\def\sectionrunning#1{\setbox\sectionbox\hbox{\eightpoint #1}
  \immediate\write1{\string\def \string\REF
      \romannumeral\numbersection \string{%
      \noexpand#1 \string\dotfill \space \number\pageno \string}}}
\def\section#1{%
  \par\vskip0.666cm\penalty -100
  \vbox{\baselineskip=14.4pt\noindent{{\twelvepointbf #1}}}
  \vskip2pt
  \penalty 500
  \advance\numbersection by 1
  \sectionrunning{#1}}

\def\subsection#1{%
  \par\vskip0.5cm\penalty -100
  \vbox{\noindent{{\tenpointbf #1}}}
  \vskip1pt
  \penalty 500}

\newcount\numberindex \numberindex=0
\def\index#1#2{%
  \advance\numberindex by 1
  \immediate\write1{\string\def \string\IND #1%
     \romannumeral\numberindex \string{%
     \noexpand#2 \string\dotfill \space \string\S \number\numbersection,
     p.\string\ \space\number\pageno \string}}}

\newdimen\itemindent \itemindent=\parindent

\def\item#1{\par\noindent\hangindent\itemindent%
            \rlap{#1}\kern\itemindent\ignorespaces}
\def\itemitem#1{\par\noindent\hangindent2\itemindent%
            \kern\itemindent\rlap{#1}\kern\itemindent\ignorespaces}
\def\itemitemitem#1{\par\noindent\hangindent3\itemindent%
            \kern2\itemindent\rlap{#1}\kern\itemindent\ignorespaces}

\long\def\claim#1|#2\endclaim{\par\vskip 5pt\noindent
{\tenpointbf #1.}\ {\it #2}\par\vskip 5pt}

\def\proof{\noindent{\it Proof}}

\def\today{\ifcase\month\or
January\or February\or March\or April\or May\or June\or July\or August\or
September\or October\or November\or December\fi \space\number\day,
\number\year}

\catcode`\@=11
\newcount\@tempcnta \newcount\@tempcntb
\def\timeofday{{%
\@tempcnta=\time \divide\@tempcnta by 60 \@tempcntb=\@tempcnta
\multiply\@tempcntb by -60 \advance\@tempcntb by \time
\ifnum\@tempcntb > 9 \number\@tempcnta:\number\@tempcntb
  \else\number\@tempcnta:0\number\@tempcntb\fi}}
\catcode`\@=12

\def\bibitem#1&#2&#3&#4&%
{\hangindent=1.8cm\hangafter=1
\noindent\rlap{\hbox{\eightpointbf #1}}\kern1.8cm{\rm #2}{\it #3}{\rm #4.}}


\def\bC{{\Bbb C}}

\def\bN{{\Bbb N}}
\def\bP{{\Bbb P}}
\def\bQ{{\Bbb Q}}
\def\bR{{\Bbb R}}

\def\bZ{{\Bbb Z}}

\def\cE{{\Cal E}}
\def\cF{{\Cal F}}

\def\cK{{\Cal K}}
\def\cI{{\Cal I}}
\def\cJ{{\Cal J}}
\def\cL{{\Cal L}}
\def\cO{{\Cal O}}
\def\cR{{\Cal R}}
\def\cS{{\Cal S}}


\def\bu{{\scriptstyle\bullet}}

\def\square{{\hfill \hbox{
\vrule height 1.453ex  width 0.093ex  depth 0ex
\vrule height 1.5ex  width 1.3ex  depth -1.407ex\kern-0.1ex
\vrule height 1.453ex  width 0.093ex  depth 0ex\kern-1.35ex
\vrule height 0.093ex  width 1.3ex  depth 0ex}}}
\def\qed{\kern10pt$\square$}
\def\hexnbr#1{\ifnum#1<10 \number#1\else
 \ifnum#1=10 A\else\ifnum#1=11 B\else\ifnum#1=12 C\else
 \ifnum#1=13 D\else\ifnum#1=14 E\else\ifnum#1=15 F\fi\fi\fi\fi\fi\fi\fi}
\def\msatype{\hexnbr\msafam}
\def\msbtype{\hexnbr\msbfam}
\mathchardef\restriction="3\msatype16   
\mathchardef\boxsquare="3\msatype03
\mathchardef\preccurlyeq="3\msatype34
\mathchardef\compact="3\msatype62
\mathchardef\smallsetminus="2\msbtype72   \let\ssm\smallsetminus
\mathchardef\subsetneq="3\msbtype28
\mathchardef\supsetneq="3\msbtype29
\mathchardef\leqslant="3\msatype36   \let\le\leqslant
\mathchardef\geqslant="3\msatype3E   \let\ge\geqslant
\mathchardef\stimes="2\msatype02
\mathchardef\ltimes="2\msbtype6E
\mathchardef\rtimes="2\msbtype6F

\def\dbar{\overline\partial}
\def\ddbar{\partial\overline\partial}


\let\la=\longrightarrow
\def\ul#1{$\underline{\hbox{#1}}$}
\let\wt=\widetilde
\let\wh=\widehat
\let\text=\hbox
\def\buildo#1^#2{\mathop{#1}\limits^{#2}}
\def\buildu#1_#2{\mathop{#1}\limits_{#2}}
\def\ort{\mathop{\hbox{\kern1pt\vrule width4.0pt height0.4pt depth0pt
                \vrule width0.4pt height6.0pt depth0pt\kern3.5pt}}}
\let\lra\longrightarrow
\def\vlra{\mathrel{\smash-}\joinrel\mathrel{\smash-}\joinrel%
\kern-2pt\longrightarrow}

\def\Re{\mathop{\rm Re}\nolimits}

\def\Alb{\mathop{\rm Alb}\nolimits}
\def\Pic{\mathop{\rm Pic}\nolimits}

\def\eff{\mathop{\rm eff}\nolimits}
\def\nef{\mathop{\rm nef}\nolimits}
\def\ample{\mathop{\rm ample}\nolimits}
\def\rk{\mathop{\rm rk}\nolimits}
\def\codim{\mathop{\rm codim}\nolimits}
\def\Id{\mathop{\rm Id}\nolimits}
\def\Sing{\mathop{\rm Sing}\nolimits}

\def\Todd{\mathop{\rm Todd}\nolimits}


\def\cH{{\cal H}}

\long\def\InsertFig#1 #2 #3 #4\EndFig{\par
\hbox{\hskip #1mm$\vbox to#2mm{\vfil\special{"
(/home/demailly/psinputs/grlib.ps) run
#3}}#4$}}
\long\def\LabelTeX#1 #2 #3\ELTX{\rlap{\kern#1mm\raise#2mm\hbox{#3}}}


\def\REF{\S 0. Introduction \dotfill 1}

\let\sectref=\SECTREF


\title{Pseudo-effective line bundles}
\title{on compact K\"ahler manifolds}
\titlerunning{Pseudo-effective line bundles on compact K\"ahler manifolds}
\vskip10pt

{\noindent\hangindent0.6cm\hangafter-1
{\bf Jean-Pierre Demailly${}^\star$, Thomas Peternell${}^{\star\star}$,
Michael Schneider${}^{\star\star}$

}}

{\noindent\hangindent0.6cm\hangafter-4{\sl
\llap{${}^\star~$}Universit\'e de Grenoble I, BP 74\hfill
${}^{\star\star}~$Universit\"at Bayreuth\kern0.6cm\break
Institut Fourier, UMR 5582 du CNRS
\hfill Mathematisches Institut\kern0.6cm\break
38402 Saint-Martin d'H\`eres, France\hfill D-95440 Bayreuth,
Deutschland\kern0.6cm

}}
\vskip20pt

{\eightpoint
\noindent{\bf Abstract.} The goal of this work is to pursue the study of
pseudo-effective line bundles and vector bundles. Our first result is a
generalization of the Hard Lefschetz theorem for cohomology with
values in a pseudo-effective line bundle. The Lefschetz map is shown
to be surjective when (and, in general, only when) the pseudo-effective
line bundle is twisted by its multiplier ideal sheaf. This result has
several geometric applications, e.g.\ to the study of compact K\"ahler
manifolds with pseudo-effective canonical or anti-canonical line
bundles. Another concern is to understand pseudo-effectivity in more
algebraic terms.  In this direction, we introduce the concept of an
``almost'' nef line bundle, and mean by this that the degree of the
bundle is nonnegative on sufficiently generic curves. It can be shown
that pseudo-effective line bundles are almost nef, and our hope is that
the converse also holds true. This can be checked in some
cases, e.g.\ for the canonical bundle of a projective $3$-fold. From this,
we derive some geometric properties of the Albanese map of compact
K\"ahler $3$-folds.\par}  \vskip14pt

\noindent{\bf Contents}
\vskip5pt

{\eightpoint
\sectref{0}
\sectref{1}
\sectref{2}
\sectref{3}
\sectref{4}
\sectref{5}
\sectref{6}
\sectref{7}}
\vskip8pt

\section{\S0. Introduction}

A line bundle $L$ on a projective manifold $X$ is {\it
  pseudo-effective} if $c_1(L) $ is in the closed cone in
$H^{1,1}_\bR(X)$ generated by the effective divisors. If $X$ is
only supposed to be K\"ahler, this definition is no longer very
meaningful, instead we require that $c_1(L)$ is in the closure of
the cone generated by
the classes of closed positive $(1,1)$-currents. In case $X$ is
projective this is equivalent to the previous definition.
Pseudo-effective line bundles on K\"ahler manifolds were first
introduced in [De90]. The aim of this paper is to study
pseudo-effective line bundles in general (as well as the concept
of pseudo-effective vector bundle) and in particular varieties
whose canonical or anticanonical
bundles are pseudo-effective. Pseudo-effective line bundles can also
be characterized in a differential-geometric way: they carry singular
hermitian metrics $h$, locally of the form $e^{-2\varphi}$ with
$\varphi$ integrable, such that the curvature current
$$\Theta_h(L) = i \partial {\overline {\partial}}\varphi$$
is positive. In this context the multiplier ideal sheaf $\cI(h) $ plays an
important role; by definition it is the ideal sheaf of local
holomorphic functions $f$ such that $\vert f \vert^2e^{-2\varphi}$ is
locally integrable.  If $h$ is a smooth metric (with semi-positive
curvature), then $\cI(h) = \cO_X$, but the converse is not true.
Our first main result in $\S$~2 is the following hard Lefschetz theorem

\claim Theorem|Let $(L,h)$ be a pseudo-effective line bundle on a
compact K\"ahler manifold $(X,\omega)$ of dimension~$n$, let
$\Theta_h(L)\ge 0$ be its
curvature current and $\cI(h)$ the associated multiplier ideal sheaf.
Then, for every nonnegative integer $q$, the wedge multiplication
operator $\omega^q\wedge\bu$ induces a surjective morphism
$$
H^0(X,\Omega_X^{n-q}\otimes L\otimes\cI(h))\buildo\vlra^{\textstyle
\omega^q\wedge\bu} H^q(X,\Omega_X^n\otimes L\otimes\cI(h)).
$$
\endclaim

The classical hard Lefschetz theorem is the case when $L$ is trivial or
unitary flat;
then $L$ has a metric $h$ of zero curvature, whence $\cI(h) =
\cO_X$. One might ask whether it is possible to
omit the multiplier ideal sheaf when $L$ is nef, i.e.\ when $c_1(L)$
is contained in the closure of the K\"ahler cone. This is however not
the case, as shown by some example (see $\S$~2.5). Therefore the above
Lefschetz theorem provides a tool to distinguish between nef and hermitian
semi-positive line bundles at the cohomological level.

\medskip We then give two applications of the hard Lefschetz theorem.
The first is concerned with compact K\"ahler manifolds $X$ whose
anticanonical bundle $-K_X$ is pseudo-effective, carrying a singular
metric with semi-positive curvature and $\cI(h) = \cO_X.$ Then we show
that the Albanese map of $X$ is a surjective submersion. We will come
back to this type of problems later. The second application deals with
K\"ahler manifolds with $K_X$ pseudo-effective and should be
considered as a contribution to Mori theory on K\"ahler manifolds.  We
show that if $K_X$ has a singular metric $h$ whose singularities are
not ``too bad" (e.g.\ if $\cI(h) = \cO_X))$, then either
$\chi(X,\cO_X) = 0$ (which provides a non-zero holomorphic form of odd
degree) or
$$
H^0(X,\Omega_X^q \otimes \cO_X(mK_X)) \ne 0
$$
for a fixed number $q$ and infinitely many $m.$ One might hope that
this last condition enforces $\kappa (X) \geq 0;$ at least we are able
to show that it implies the existence of a non-constant meromorphic
function on $X.$ Using results of threefold classification, we are
then able to show $\kappa (X) \geq 0$ for a compact K\"ahler threefold
with $K_X$ pseudo-effective having a metric of semi-positive curvature
with ``mild" singularities. In particular this holds if $K_X$ is
hermitian semi-positive. Of course, the algebraic case (``Abundance
Conjecture") is known since some time by deep results of Mori theory.
We also prove, however not as an application of the hard Lefschetz
theorem and therefore postponed to $\S$~5, that a compact K\"ahler
threefold (isolated singularities are allowed but $X$ must be
$\bQ$-factorial) with $K_X$ pseudo-effective but not nef admits a curve
$C$ with $K_X \cdot C < 0.$ In case $X$ is smooth this implies the
existence of a Mori contraction. Of course this is only new in the
non-algebraic setting.

\medskip We next address (in $\S$~3) the question whether
pseudo-effective line bundles can be characterized
in more algebraic terms in case the underlying manifold is
projective. We say that a line bundle $L$ is {\it almost nef},
if there is a family $A_i \subset X$,
$i\in\bN$, of proper algebraic subvarieties such that $L\cdot C \ge 0$
for all irreducible curves $C \not\subset \bigcup_{i} A_i$.
The Zariski closure of the union of all curves $C$ with $L\cdot C<0$
will be called the {\it non-nef locus} of $L$.
In this setting, pseudo-effective line bundles $L$ turn out
to be almost nef, but the converse seems to be a very hard problem (if
at all true). The equivalence between pseudo-effectivity and almost nefness
is however always true on surfaces; and, by using Mori theory, it is true
for $L = K_X$ on every threefold.

\medskip In $\S$~4 we study compact K\"ahler manifolds $X$
with $-K_X$ pseudo-effective resp. almost nef. First we study morphism
$\varphi: X \longrightarrow Y$ and restrict to projective varieties
$X$ and $Y$. In general $-K_Y$ will not be pseudo-effective resp.
almost nef; the reason is that the non-nef locus of $-K_X$ might project
onto $Y$. Ruling this out, we obtain

\claim Theorem|Let $X$ and $Y$ be normal projective $\bQ$-Gorenstein
varieties. Let $\varphi: X \to Y$ be a surjective map with connected
fibers.
\smallskip
\item{\rm(a)} Suppose that $X$ and $Y$ are smooth, that
  $\varphi$ is a submersion, that $-K_X$ is pseudo-effective and that
  the zero locus of the multiplier ideal of a minimal metric on $-K_X$ does
not project onto $Y.$ Then
  $-K_Y$ is pseudo-effective.
\smallskip
\item{\rm(b)}  Let $-K_X$
be almost nef with non-nef locus not projecting onto $Y$.
Then $-K_Y$ is generically nef.
\vskip0pt
\endclaim

We say that a  $(\bQ$-) line bundle $L$ on a normal
$n$-dimensional projective variety $X$ is {\it generically nef} if
$$L\cdot H_1 \cdots H_{n-1} \ge 0$$
for all ample divisors $H_i$ on
$X$. This is a much weaker notion than almost nefness.

\claim  Corollary|Let $X$ be a normal projective
$\bQ$-Gorenstein variety. Assume $-K_X$ almost
nef with non-nef locus $B$.
\smallskip
\item{\rm (a)} If $\varphi:X \to Y$ is a surjective morphism to a
  normal projective $\bQ$-Gorenstein variety $Y$ with
  $\varphi(B) \not= Y$, then $\kappa (Y) \le 0.$
\smallskip
\item{\rm (b)} The Albanese map $\alpha:X \to \Alb(X)$ is surjective,
  if $\alpha (B) \not= \alpha (X).$ \vskip0pt
\endclaim

We shall call $-K_X$ {\it properly pseudo-effective
$($properly almost nef$)$} with respect to $\alpha$, if
it is pseudo-effective (almost nef), with
$\alpha (B) \not= \alpha (X)$. Concerning the fundamental group,
we show

\claim Theorem|Let $X$ be a terminal 3-fold with $-K_X$ properly
pseudo-effective $($properly almost nef$)$ with respect to the Albanese
map. Then $\pi_1(X)$ is almost abelian.
\endclaim

If $-K_X$ is hermitian semi-positive, then $\pi_1(X)$ is almost abelian
by [DPS96b] (this is true in any dimension); in case $-K_X$ nef,
we know by (Paun [Pa96]) that $\pi_1(X)$ has polynomial growth.

\medskip If $X$ is merely supposed to be K\"ahler, we have
to restrict to the case when $-K_X$ is nef.
Here, using the algebraic case settled in [PS97], our main result is

\claim Theorem|Let $X$ be a compact K\"ahler 3-fold with $-K_X$
nef. Then the Albanese map is a surjective submersion. Moreover
$\pi_1(X)$ is almost abelian.
\endclaim

The last section is concerned with pseudo-effective and almost nef
vector bundles.  Given a projective manifold $X$ and a vector bundle
$E$ over $X$, we say that $E$ is {\it pseudo-effective}, if
$\cO_{\bP(E)}(1)$ is pseudo-effective and the union of all curves $C$
with $\cO(1) \cdot C < 0$ $($i.e.\ the non-nef locus of the almost nef
line bundle $\cO(1))$ is contained in a union of subvarieties which does not
project onto $X$.  The definition of almost nefness is analogous to
the rank 1 case.  Basic properties of pseudo-effective and almost nef
bundles are collected in

\claim Theorem|Let $X$ be a projective manifold and $E$ a vector
bundle on $X.$
\smallskip
\item{\rm (a)} If $E$ is pseudo-effective (almost nef) and $\Gamma^{a}$
  is any tensor representation, then $\Gamma^{a}E$ is again
  pseudo-effective $($almost nef$)$. In particular $S^mE$ and $\Lambda^qE$
  are pseudo-effective $($almost nef$)$.  \smallskip
\item{\rm (b)} If $E$ is almost nef and if $s \in H^0(E^*)$ is a non-zero
section, then $s$ has no
zeroes at all.
\smallskip
\item{\rm (c)}
  If either $E$ is pseudo-effective or almost nef with non-nef locus
  $S$ having codimension at least $2,$ and if $\det E^*$ is almost
  nef, then $E$ is numerically flat i.e.\ both $E$ and $E^*$ are nef,
  and then $E$ has a filtration by hermitian flat bundles.
\vskip0pt
\endclaim

We then discuss projective manifolds $X$ with pseudo-effective or
almost nef tangent bundles $T_X$. Important examples are provided by
almost homogeneous spaces, i.e.\ the automorphism group acts with an
open orbit. It is easily seen that necessarily $\kappa (X) \leq 0.$
Moreover if $\kappa (X) = 0,$ then $K_X \equiv 0.$ Now in case $K_X
\equiv 0,$ and if $T_X$ is pseudo-effective or if $T_X$ is almost nef
with non-nef locus of codimension at least 2, then $X$ is abelian
after possibly finite \'etale cover. As a consequence, a Calabi-Yau
manifold can never have a pseudo-effective tangent or cotangent bundle
and the union of the curves $C$ with $T_X \vert C$ is not nef (resp.
$\Omega^1_X \vert C$ is not nef) is not contained in a countable union
of analytic sets of codimension at least 2.
\medskip

In general the Albanese map of $X$ are surjective submersion and we give
a precise description of the Albanese map in case $\dim X = 3.$ In
case neither $X$ nor any finite \'etale cover of $X$ has a holomorphic
1-form, $X$ is expected to be simply connected. This is true if $X$ is
almost homogeneous, if $\dim X \leq 3$ or if $T_X$ is nef ([DPS94]).

\medskip The research on this paper started in fall 1996
with the important participation of Michael Schneider. After his
tragic death in august 1997, the paper was finished -- after some delay
-- by the two first-named authors who therefore carry full scientific
responsibility. We would like to thank the referee for his very careful
reading of the manuscript.

\section{\S1. Nef and pseudo-effective line bundles}

We first recall (in the K\"ahler context) the basic concepts of
numerical effectivity and pseudo-effectivity. The proofs as well as
more details can be found in [De90, De92, DPS94, DPS96a]. Given a
holomorphic line bundle on a complex manifold $X$ and a hermitian
metric $h$ on $L$, we denote by $\Theta_h(L)=i D_h^2$ the curvature
of the Chern connection $D_h$ associated with~$h$. This is a real
$(1,1)$-form, which can be expressed as $\Theta_h(L)=-i\ddbar\log h$
in coordinates. The first Chern class of $L$ is represented by
$\{{1\over 2\pi}\Theta_h(L)\}$ in $H^{1,1}(X)\subset H^2(X,\bR)$.

\claim 1.1.\ Definition|Let $X$ be a compact K\"ahler manifold. A line
bundle $L$ on $X$ is said to be
\smallskip
\item{\rm a)}
{\it pseudo-effective} if $c_1(L)$ is in the closed cone of $H^{1,1}_\bR(X)$
generated by classes of $d$-closed positive $(1,1)$-currents.
\smallskip
\item{\rm b)} {\it nef} $($numerically effective$)$ if $c_1(L)$ is in the
closure of the K\"ahler cone, i.e.\ the closed cone generated by smooth
non-negative $d$-closed $(1,1)$-forms.
\vskip0pt
\endclaim

\noindent
It is clear from the above definition that every nef line bundle is
pseudo-effective (but the converse is in general {\it not true}).
The names of these concepts stem from the following ``more concrete''
characterisation  in case $X$ is projective.

\claim 1.2.\ Proposition|Let $X$ be a projective manifold and $L$ a
line bundle on $X$. Then
\smallskip
\item{\rm a)} $L$ is pseudo-effective if and only if $c_1(L)$ is
in the closure $\bar{K}_{\eff} (X)$ of the cone generated by the
effective divisors $($modulo numerical equivalence$)$ on $X$.
\smallskip
\item{\rm b)} $L$ is nef if and only if the degree $L\cdot C$ is
non-negative for every effective curve $C\subset X$, or equivalently,
if $c_1(L)$ is in the closure $K_{\nef}(X)=\bar{K}_{\ample}$ of the
cone of ample divisors.
\vskip0pt
\endclaim

Assume that $L$ is pseudo-effective, and let $T$ be a closed positive
$(1,1)$-current such that $c_1(L) = [T].$ Choose a smooth hermitian metric
$h_\infty$ on $L$. Let $\alpha={1\over 2\pi}\Theta_{h_\infty}(L)$ denote
its curvature. Since $\{T\}=\{\alpha\}$, we can write
$T = \alpha + {i \over \pi}\ddbar\psi$ for some locally integrable
function $\psi$. Then $h = h_\infty \exp(-2\psi)$ is a singular
metric on $L$, and
$$
{1\over 2\pi}\Theta_h(L)=\alpha + {i \over \pi}\ddbar\psi = T.
$$
Hence, $L$ is pseudo-effective if and only if
there exists a singular hermitian metric $h$ on $L$ such that
its curvature current $\Theta_h(L)=-i\ddbar\log h$ is positive.

When $X$ is projective, H\"ormander's $L^2$ estimates show that
$L$ is pseudo-effective if and only if there exists an ample divisor
$A$ such that
$$H^0(X, \cO_X(mL+A)) \not= 0$$
for $m \gg 0$ (see [De90]). Whence 1.2~a). We now come to the very
important concept of multiplier ideal sheaf.

\claim 1.3.\ Multiplier ideal sheaves|Let $\varphi$ be a psh
$($plurisubharmonic$)$ function on an open subset $\Omega\subset \bC^n$.
To $\varphi$ is
associated the ideal subsheaf $\cI(\varphi)\subset\cO_\Omega$ of germs of
holomorphic functions $f\in\cO_{\Omega,x}$ such that $|f|^2e^{-2\varphi}$
is integrable with respect to the Lebesgue measure in some local
coordinates near~$x$.
\endclaim

A basic result of Nadel [Nad89] shows that the sheaf $\cI(\varphi)$ is
a coherent sheaf of ideals over $\Omega$, generated by its global
$L^2$ sections over $\Omega$, provided $\Omega$ is e.g.\ bounded and
pseudoconvex (this comes from standard $L^2$ estimates combined with
the Krull lemma). If $(L,h)$ is a pseudo-effective line bundle
with $\Theta_h(L)\ge 0$, then the local weight functions $\varphi$
of $h$ are plurisubharmonic and we simply denote $\cI(h):=\cI(\varphi)$.

\claim 1.4.\ Definition|Let $L$ be a pseudo-effective line bundle
on a compact complex manifold~$X$. Consider two hermitian
metrics $h_1$, $h_2$ on $L$ with curvature $\Theta_{h_j}(L)\ge 0$
in the sense of currents.
\smallskip
\item{\rm a)} We write $h_1\preccurlyeq h_2$, and say that $h_1$ is less
singular than $h_2$, if there exists a constant $C>0$ such that
$h_1\le C h_2$.
\smallskip
\item{\rm b)} We write $h_1\sim h_2$, and say that $h_1$, $h_2$ are
singularity equivalent, if there exists a constant
$C>0$ such that $C^{-1}h_2\le h_1\le C h_2$.
\vskip0pt
\endclaim

\noindent
Of course $h_1\preccurlyeq h_2$ if and only if the local associated weights
in suitable trivializations satisfy $\varphi_2\le\varphi_1+C$. This implies
in particular that the lelong numbers satisfy
$\nu(\varphi_1,x)\le\nu(\varphi_2,x)$ at every point. The above definition
is motivated by the following observation.

\claim 1.5.\ Theorem|For every pseudo-effective line bundle $L$ over
a compact complex manifold $X$, there exists up to equivalence of
singularities a unique class of hermitian metrics $h$ with minimal
singularities such that $\Theta_h(L)\ge 0$.
\endclaim

\proof. The proof is almost trivial. We fix once for all a smooth
metric $h_\infty$ (whose curvature is of random sign and signature),
and we write singular metrics of $L$ under the form $h=h_\infty e^{-2\psi}$.
The condition $\Theta_h(L)\ge 0$ is equivalent to
${i\over\pi}\ddbar\psi\ge -u$ where $u=\Theta_{h_\infty}(L)$. This
condition implies that $\psi$ is plurisubharmonic up to the addition
of the weight $\varphi_\infty$ of $h_\infty$, and therefore locally
bounded from above.  Since we are concerned with metrics only up to
equivalence of singularities, it is always possible to adjust
$\psi$ by a constant in such a way that $\sup_X\psi=0$. We now set
$$
h_{\min}=h_\infty e^{2\psi_{\min}},\qquad \psi_{\min}(x)=\sup_\psi\psi(x)
$$
where the supremum is extended to all functions $\psi$ such that
$\sup_X\psi=0$ and\break ${i\over\pi}\ddbar\psi\ge -u$. By standard
results on plurisubharmonic functions (see Lelong [Lel69]),
$\psi_{\min}$ still satisfies ${i\over\pi}\ddbar\psi_{\min}\ge -u$
(i.e.\ the weight $\varphi_\infty+\psi_{\min}$ of $h_{\min}$ is
plurisubharmonic), and $h_{\min}$ is obviously the metric with
minimal singularities that we were looking for. \qed

Now, given a section $\sigma\in H^0(X,mL)$, the expression $h(\xi)=
|\xi^m/\sigma(x)|^{2/m}$ defines a singular metric on $L$,
which therefore necessarily has at least as much singularity as
$h_{\min}$ as,
i.e.\ ${1\over m}\log|\sigma|^2\le \varphi_{\min}+C$ locally. In particular,
$|\sigma|^2e^{-2m\varphi_{\min}}$ is locally bounded, hence
$\sigma\in H^0(X,mL\otimes\cI(h_{\min}^{\otimes m}))$. For all
$m>0$ we therefore get an isomorphism
$$
H^0(X,mL\otimes\cI(h_{\min}^{\otimes m}))\buildo\vlra^\simeq H^0(X,mL).
$$
By the well-known properties of Lelong numbers (see Skoda [Sk72]), the
union of all zero varieties of the ideals $\cI(h_{\min}^{\otimes m})$
is equal to the Lelong sublevel set
$$
E_+(h_{\min})=\big\{x\in X\,;~\nu(\varphi_{\min},x)>0\big\}.
\leqno(1.6)
$$
We will call this set the {\it virtual base locus} of $L$. It is
always contained in the ``algebraic'' base locus
$$B_{\Vert L\Vert}=\bigcap_{m>0}B_{|mL|},\qquad
B_{|mL|}=\bigcap_{\sigma\in H^0(X,mL)}\sigma^{-1}(0),
$$
but there may be a strict inclusion. This is the case for instance if
$L\in\Pic^0(X)$ is such that all positive multiples $mL$ have no
nonzero sections; in that case $E_+(h_{\min})=\emptyset$ but
$\bigcap_{m>0}B_{|mL|}=X$. Another general situation where
$E_+(h_{\min})$ and $B_{\Vert L\Vert}$ can differ is given by the
following result.

\claim 1.7.\ Proposition|Let $L$ be a big nef line bundle. Then $h_{\min}$
has zero Lelong numbers everywhere, i.e.\ $E_+(h_{\min})=\emptyset$.
\endclaim

\proof. Recall that $L$ is big if its Kodaira-Iitaka dimension
$\kappa(L)$ is equal to $n=\dim X$. In that case, it is well known that
one can write $m_0L=A+E$ with $A$ ample and $E$ effective, for $m_0$
sufficiently large. Then $mL=((m-m_0)L+A) +E$ is the sum of an
ample divisor $A_m=(m-m_0)L+A$ plus a (fixed) effective divisor,
so that there is a hermitian metric $h_m$ on $L$ for which
$\Theta_{h_m}(L)={1\over m}\Theta(A_m)+{1\over m}[E]$,
with a suitable smooth positive form $\Theta(A_m)$. This shows that
the Lelong numbers of the weight of $h_m$ are $O(1/m)$, hence in
the limit those of $h_{\min}$ are zero. \qed

If $h$ is a singular hermitian metric such that $\Theta_h(L)\ge 0$ and
$$
H^0(X,mL\otimes\cI(h^{\otimes m}))\simeq H^0(X,mL)\qquad
\hbox{for all $m\ge 0$},
\leqno(1.8)
$$
we say that $h$ is an {\it analytic Zariski decomposition} of~$L$.
We have just seen that such a decomposition always exists and
that $h=h_{\min}$ is a solution. The concept of analytic Zariski decomposition
is motivated by its algebraic counterpart (the existence of which generally
fails)$\,$: one says that $L$ admits an {\it algebraic Zariski decomposition}
if there exists a modification $\mu:\tilde{X}\to X$ and an integer $m_0$
with $m_0\tilde{L}\simeq\cO(E+D)$, where $\tilde{L}=\mu^\star L$,
$E$ is an effective divisor and $D$ a nef divisor on $\tilde{X}$ such that
$$H^0(\tilde{X},kD)=H^0(\tilde{X},k(D+E))\simeq H^0(X,km_0L)\qquad
\hbox{for all $k\ge 0$}.\leqno(1.9)
$$
If $\cO(*D)$ is generated by sections, there is a smooth metric with
semi-positive curvature on $\cO(D)$, and this metric induces a
singular hermitian metric $\tilde{h}$ on $\tilde{L}=\mu^\star(L)$ of
curvature current ${1\over m_0}(\Theta(\cO(D))+ [E])$. Its poles are
defined by the effective $\bQ$-divisor ${1\over m_0}E$.  For this
metric, we of course have $\cI(\tilde{h}^{\otimes km_0}) =\cO(-kE)$,
hence assumption (1.9) can be rewritten
$$H^0\big(\tilde{X},km_0\tilde{L}\otimes\cI(\tilde{h}^{\otimes km_0})\big)=
H^0(\tilde{X},km_0\tilde{L})\qquad
\hbox{for all $k\ge 0$}.\leqno(1.10)
$$
When we take the direct image, we find a hermitian metric $h$ on $L$
with curvature current $\Theta_h(L)=\mu_\star\Theta_{\tilde h}(
\tilde{L})\ge 0$ and
$$\cO_X\supset \cI(h^{\otimes km_0})=
\mu_\star\big(K_{\tilde{X}/X}\otimes\cI(\tilde{h}^{\otimes km_0})\big)
\supset \mu_\star\big(\cI(\tilde{h}^{\otimes km_0})\big),
\leqno(1.11)
$$
thus (1.8) holds true at least when $m$ is a multiple of $m_0$.\qed

\section{\S2. Hard Lefschetz theorem with multiplier
ideal sheaves}

\subsection{\S2.1. Main statements}

The goal of this section is to prove the following surjectivity theorem,
which can be seen as an extension of the hard Lefschetz theorem.

\claim 2.1.1.\ Theorem|Let $(L,h)$ be a pseudo-effective line bundle
on a compact K\"ahler manifold $(X,\omega)$ of dimension~$n$, let
$\Theta_h(L)\ge 0$ be its curvature current and $\cI(h)$ the
associated multiplier ideal sheaf.  Then, for every nonnegative integer
$q$, the wedge multiplication operator $\omega^q\wedge\bu$ induces a
surjective morphism
$$
\Phi^q_{\omega,h}:
H^0(X,\Omega_X^{n-q}\otimes L\otimes\cI(h))\lra
H^q(X,\Omega_X^n\otimes L\otimes\cI(h)).
$$
\endclaim

\noindent
The special case when $L$ is nef is due to Takegoshi [Ta97]. An even
more special case is when $L$ is semi-positive, i.e.\ possesses a
smooth metric with semi-positive curvature.
In that case the multiplier ideal sheaf $\cI(h)$ coincides with $\cO_X$
and we get the following consequence already observed by Mourougane
[Mou99].

\claim 2.1.2.\ Corollary|Let $(L,h)$ be a semi-positive line bundle
on a compact K\"ahler mani\-fold $(X,\omega)$ of dimension~$n$. Then, the
wedge multiplication operator $\omega^q\wedge\bu$ induces a surjective
morphism
$$
\Phi^q_\omega:H^0(X,\Omega_X^{n-q}\otimes L)\lra
H^q(X,\Omega_X^n\otimes L).
$$
\endclaim

The proof of Theorem 2.1.1 is based on the Bochner formula, combined
with a use of harmonic forms with values in the hermitian line bundle
$(L,h)$. The method can be applied only after $h$ has been made smooth
at least in the complement of an analytic set. However, we have to
accept singularities even in the regularized metrics because only a
very small incompressible loss of positivity is acceptable in the
Bochner estimate (by the results of [De92], singularities can be
removed, but only at the expense of a fixed, non zero, loss of
positivity). Also, we need the multiplier ideal sheaves to be preserved
by the smoothing process. This is possible thanks to a suitable
``equisingular'' regularization process.

\subsection{\S2.2. Equisingular approximations of quasi plurisubharmonic
functions}

A quasi-plurisubharmonic (quasi-psh) function is by definition a
function $\varphi$ which is locally equal to the sum of a psh function
and of a smooth function, or equivalently, a locally integrable function
$\varphi$ such that $i\ddbar\varphi$ is locally bounded below by
$-C\omega$ where $\omega$ is a hermitian metric and $C$ a constant.
We say that $\varphi$ has logarithmic poles if $\varphi$ is locally
bounded outside an analytic set $A$ and has singularities of the form
$$\varphi(z)=c\log\sum_k|g_k|^2+O(1)$$
with $c>0$ and $g_k$ holomorphic, on a neighborhood of every point
of~$A$. Our goal is to show the following

\claim 2.2.1.\ Theorem|Let $T=\alpha+i\ddbar\varphi$ be a closed
$(1,1)$-current on a compact hermitian manifold $(X,\omega)$, where
$\alpha$ is a smooth closed $(1,1)$-form and $\varphi$ a quasi-psh
function. Let $\gamma$ be a continuous real $(1,1)$-form such that
$T\ge\gamma$. Then one can write
\hbox{$\varphi=\lim_{\nu\to+\infty}\varphi_\nu$} where
\smallskip
\item{\rm a)} $\varphi_\nu$ is smooth in the complement $X\ssm Z_\nu$
of an analytic set $Z_\nu\subset X\,;$
\smallskip
\item{\rm b)} $(\varphi_\nu)$ is a decreasing sequence, and $Z_\nu\subset
Z_{\nu+1}$ for all~$\nu\,;$
\smallskip
\item{\rm c)} $\int_X(e^{-2\varphi}-e^{-2\varphi_\nu})dV_\omega$
is finite for every $\nu$ and converges to $0$ as $\nu\to+\infty\,;$
\smallskip
\item{\rm d)} $\cI(\varphi_\nu)=\cI(\varphi)$ for all $\nu$
$($``equisingularity''$)\,;$
\smallskip
\item{\rm e)} $T_\nu=\alpha+i\ddbar\varphi_\nu$ satisfies $T_\nu\ge
\gamma-\varepsilon_\nu\omega$, where $\lim_{\nu\to+\infty}\varepsilon_\nu=0$.
\vskip0pt
\endclaim

\claim 2.2.2.\ Remark|{\rm It would be interesting to know whether the
$\varphi_\nu$ can be taken to have logarithmic poles along $Z_\nu$.
Unfortunately, the proof given below destroys this property in the
last step. Getting it to hold true seems to be more or less
equivalent to proving the semi-continuity property
$$\lim_{\varepsilon\to 0_+}
\cI((1+\varepsilon)\varphi)=\cI(\varphi).$$
Actually, this can be checked in dimensions $1$ and $2$, but is unknown in
higher dimensions (and probably quite hard to establish).}
\endclaim

\proof\ {\it of Theorem 2.2.1}. Clearly, by replacing $T$ with
$T-\alpha$ and $\gamma$ with $\gamma-\alpha$, we may assume that
$\alpha=0$ and $T=i\ddbar\varphi\ge \gamma$. We divide the proof in
four steps.
\medskip

\noindent{\it Step 1. Approximation by quasi-psh functions with
logarithmic poles.}

\noindent
By [De92], there is a decreasing sequence $(\psi_\nu)$ of quasi-psh
functions with logarithmic poles such that $\varphi=\lim\psi_\nu$ and
$i\ddbar\psi_\nu\ge\gamma- \varepsilon_\nu\omega$. We need a little bit
more information on those functions, hence we first recall the main
techniques used for the construction of~$(\psi_\nu)$. For $\varepsilon>0$
given, fix a covering of $X$ by open balls $B_j=\{|z^{(j)}|<r_j\}$ with
coordinates $z^{(j)}=(z^{(j)}_1,\ldots, z^{(j)}_n)$, such that
$$
0\le\gamma+c_j\,i\ddbar|z^{(j)}|^2\le \varepsilon\omega
\qquad\hbox{on}\quad B_j,
\leqno(2.2.3)
$$
for some real number~$c_j$. This is possible by selecting coordinates
in which $\gamma$ is diagonalized at the center of the ball, and by
taking the radii $r_j>0$ small enough (thanks to the fact that $\gamma$
is continuous). We may assume that these coordinates come from a finite
sample of coordinates patches covering $X$, on which we perform suitable
linear coordinate changes (by invertible matrices lying in some
compact subset of the complex linear group). By taking additional balls,
we may also assume that $X=\bigcup B''_j$ where
$$
B''_j\compact B'_j\compact B_j
$$
are concentric balls $B'_j=\{|z^{(j)}|<r'_j=r_j/2\}$,
$B''_j=\{|z^{(j)}|<r''_j=r_j/4\}$. We define
$$
\psi_{\varepsilon,\nu,j}=
{1\over 2\nu}\log\sum_{k\in\bN}|f_{\nu,j,k}|^2-c_j|z^{(j)}|^2
\qquad\hbox{on}\quad B_j,
\leqno(2.2.4)
$$
where $(f_{\nu,j,k})_{k\in\bN}$ is an orthonormal basis of the Hilbert
space ${\cal H}_{\nu,j}$ of holomorphic functions on $B_j$ with finite
$L^2$ norm
$$
\Vert u\Vert^2 = \int_{B_j} |u|^2 e^{-2\nu(\varphi+c_j|z^{(j)}|^2)}
d\lambda(z^{(j)}).
$$
(The dependence of $\psi_{\varepsilon,\nu,j}$ on $\varepsilon$ is through
the choice of the open covering~$(B_j)$). Observe that the choice of
$c_j$ in (2.2.3) guarantees that $\varphi+c_j|z^{(j)}|^2$ is
plurisubharmonic on $B_j$, and notice also that
$$
\sum_{k\in\bN}|f_{\nu,j,k}(z)|^2=\sup_{f\in\cH_{\nu,j},\,\Vert f\Vert\le 1}
|f(z)|^2
\leqno(2.2.5)
$$
is the square of the norm of the continuous linear form $\cH_{\nu,j}\to\bC$,
$f\mapsto f(z)$. We claim that there exist
constants $C_i$, $i=1,2,\ldots$ depending only on $X$ and $\gamma$
(thus independent of $\varepsilon$ and $\nu$), such that the following
uniform estimates hold:
$$
\leqalignno{\qquad\qquad
&i\ddbar\psi_{\varepsilon,\nu,j}\ge -c_j\,i\ddbar|z^{(j)}|^2\ge\gamma-
\varepsilon\omega\qquad\hbox{on}~~B'_j~~(B'_j\compact B_j),
&(2.2.6)\cr
&\varphi(z)\le\psi_{\varepsilon,\nu,j}(z)\le\sup_{|\zeta-z|\le r}\varphi(\zeta)
+{n\over\nu}\log{C_1\over r}+C_2r^2
\qquad\hbox{$\forall z\in B'_j$,~ $r<r_j-r'_j$,}
&(2.2.7)\cr
&|\psi_{\varepsilon,\nu,j}-\psi_{\varepsilon,\nu,k}|\le {C_3\over\nu}
+C_4\varepsilon\big(\min(r_j,r_k)\big)^2\qquad
\hbox{on}~~B'_j\cap B'_k.
&(2.2.8)\cr}
$$
Actually, the Hessian estimate (2.2.6) is obvious from (2.2.3) and (2.2.4).
As in the proof of ([De92], Prop.~3.1), (2.2.7) results from the
Ohsawa-Takegoshi $L^2$ extension theorem (left hand inequality)
and from the mean value inequality (right hand inequa\-lity).
Finally, as in ([De92], Lemma~3.6 and Lemma~4.6), (2.2.8) is a
consequence of H\"ormander's $L^2$ estimates. We briefly sketch the idea.
Assume that the balls $B_j$ are small enough, so that the
coordinates $z^{(j)}$ are still defined on a neighborhood of all
balls $\overline B_k$ which intersect~$B_j$ (these coordinates can be
taken to be linear transforms of coordinates belonging to a fixed
finite set of coordinate patches covering~$X$, selected once for all).
Fix a point
$z_0\in B'_j\cap B'_k$. By (2.2.4) and (2.2.5), we have
$$
\psi_{\varepsilon,\nu,j}(z_0)={1\over \nu}\log|f(z_0)|-c_j|z^{(j)}|^2
$$
for some holomorphic function $f$ on $B_j$ with $\Vert f\Vert=1$. We
consider the weight function
$$
\Phi(z)=2\nu(\varphi(z)+c_k|z^{(k)}|^2)+2n\log|z^{(k)}-z^{(k)}_0|,
$$
on both $B_j$ and $B_k$. The trouble is that we a priori have to deal with
different weights, hence a comparison of weights is needed. By the Taylor
formula applied at $z_0$, we get
$$
\Big|c_k|z^{(k)}-z^{(k)}_0|^2 - c_j|z^{(j)}-z^{(j)}_0|^2\Big|\le
C\varepsilon\big(\min(r_j,r_k)\big)^2\qquad\hbox{on $B_j\cap B_k$}
$$
[the only nonzero term of degree 2 has type $(1,1)$ and its Hessian satisfies
$$
-\varepsilon\omega\le i\ddbar(c_k|z^{(k)}|^2-c_j|z^{(j)}|^2)\le
\varepsilon\omega
$$
by (2.2.3); we may suppose $r_j\ll\varepsilon$ so that the terms of order
$3$ and more are negligible]. By writing
$|z^{(j)}|^2=|z^{(j)}-z^{(j)}_0|^2+|z^{(j)}_0|^2+
2\Re\langle z^{(j)}-z^{(j)}_0,z^{(j)}_0\rangle$, we obtain
$$\eqalign{
c_k|z^{(k)}|^2-c_j|z^{(j)}|^2
&=2c_k\Re\langle z^{(k)}-z^{(k)}_0,z^{(k)}_0\rangle
-2c_j\Re\langle z^{(j)}-z^{(j)}_0,z^{(j)}_0\rangle\cr
&\qquad{}+c_k|z^{(k)}_0|^2-c_j|z^{(j)}_0|^2
\pm C\varepsilon(\min(r_j,r_k))^2.\cr}
$$
We use a cut-off function $\theta$ equal to $1$ in a neighborhood of
$z_0$ and with support in $B_j\cap B_k$; as $z_0\in B'_j\cap B'_k$, the
function $\theta$ can be taken to have its derivatives uniformly bounded
when $z_0$ varies. We solve the equation
$\dbar u=\dbar(\theta f e^{\nu g})$ on $B_k$, where $g$ is the holomorphic
function
$$
g(z)=c_k\langle z^{(k)}-z^{(k)}_0,z^{(k)}_0\rangle
-c_j\langle z^{(j)}-z^{(j)}_0,z^{(j)}_0\rangle.
$$
Thanks to H\"ormander's $L^2$ estimates [H\"or66], the $L^2$ solution
for the weight $\Phi$ yields a holomorphic function
$f'=\theta fe^{\nu g}-u$ on $B_k$ such that $f'(z_0)=f(z_0)$ and
$$
\eqalign{
&\int_{B_k}|f'|^2e^{-2\nu(\varphi+c_k|z^{(k)}|^2)}
d\lambda(z^{(k)})
\le C'\int_{B_j\cap B_k} |f|^2 |e^{\nu g}|^2e^{-2\nu(\varphi+c_k|z^{(k)}|^2)}
d\lambda(z^{(k)})\cr
&\le
C'\exp\!\Big(\!2\nu\big(c_k|z^{(k)}_0|^2-c_j|z^{(j)}_0|^2
+C\varepsilon(\min(r_j,r_k))^2\big)\!\Big)\!
\int_{B_j}\kern-1pt|f|^2 e^{-2\nu(\varphi+c_j|z^{(j)}|^2)}
d\lambda(z^{(j)}).\cr}
$$
Let us take the supremum of ${1\over\nu}\log|f(z_0)|={1\over\nu}\log|f'(z_0)|$
over all $f$ with $\Vert f\Vert\le 1$. By the definition of
$\psi_{\varepsilon,\nu,k}$ ((2.2.4) and (2.2.5)) and the bound on
$\Vert f'\Vert$, we find
$$
\psi_{\varepsilon,\nu,k}(z_0)\le \psi_{\nu,j}(z_0)+{\log C'\over 2\nu}
+C\varepsilon(\min(r_j,r_k))^2,
$$
whence (2.2.8) by symmetry. Assume that $\nu$ is so large that
$C_3/\nu<C_4\varepsilon(\inf_j r_j)^2$.
We ``glue'' all functions $\psi_{\varepsilon,\nu,j}$
into a function $\psi_{\varepsilon,\nu}$ globally defined on~$X$,
and for this we set
$$
\psi_{\varepsilon,\nu}(z)=\sup_{j,~B'_j\ni z}
\Big(\psi_{\varepsilon,\nu,j}(z)+
12\,C_4\varepsilon(r^{\prime 2}_j-|z^{(j)}|^2)\Big)
\qquad\hbox{on}~~X.
$$
Every point of $X$ belongs to some ball $B''_k$, and for such a point
we get
$$
12\,C_4\varepsilon(r^{\prime 2}_k-|z^{(k)}|^2)\ge
12\,C_4\varepsilon(r^{\prime 2}_k-r^{\prime\prime2}_k)>2C_4r_k^2>
{C_3\over\nu}+C_4\varepsilon(\min(r_j,r_k))^2.
$$
This, together with (2.2.8), implies that in $\psi_{\varepsilon,\nu}(z)$
the supremum is never reached for indices $j$ such that
$z\in \partial B'_j$, hence $\psi_{\varepsilon,\nu}$ is well defined
and continuous, and by standard properties of upper envelopes of
(quasi)-plurisubharmonic functions we get
$$
i\ddbar\psi_{\varepsilon,\nu}\ge \gamma-C_5\varepsilon\omega
\leqno(2.2.9)
$$
for $\nu\ge\nu_0(\varepsilon)$ large enough. By inequality (2.2.7)
applied with $r=e^{-\sqrt{\nu}}$, we see that
$\lim_{\nu\to+\infty}\psi_{\varepsilon,\nu}(z) = \varphi(z)$.
At this point, the difficulty is to show that
$\psi_{\varepsilon,\nu}$ decreasing with $\nu$ -- this may not be true
formally, but we will see at Step 3 that this is essentially true.
Another difficulty is that we must simultaneously let $\varepsilon$
go to $0$, forcing us to change the covering as we want the error
to get smaller and smaller in (2.2.9).
\medskip

\noindent{\it Step 2. A comparison of integrals.}

\noindent We claim that
$$
I:=\int_X\big(e^{-2\varphi}-e^{-2\max(\varphi,{\ell\over\ell-1}
\psi_{\nu,\varepsilon})+a}\big)dV_\omega<+\infty
\leqno(2.2.10)
$$
for every $\ell\in{}]1,\nu]$ and $a\in\bR$. In fact
$$
\eqalign{
I&\le\int_{\{\varphi<{\ell\over\ell-1}\psi_{\varepsilon,\nu}+a\}}
e^{-2\varphi}dV_\omega
=\int_{\{\varphi<{\ell\over\ell-1}\psi_{\varepsilon,\nu}\}+a}
e^{2(\ell-1)\varphi-2\ell\varphi}dV_\omega\cr
&\le e^{2(\ell-1)a}\int_Xe^{2\ell(\psi_{\varepsilon,\nu}-\varphi)}dV_\omega
\le C\Big(\int_Xe^{2\nu(\psi_{\varepsilon,\nu}-\varphi)}
dV_\omega\Big)^{\ell\over\nu}\cr}
$$
by H\"older's inequality. In order to show that these integrals are finite,
it is enough, by the definition and properties of the functions
$\psi_{\varepsilon,\nu}$ and $\psi_{\varepsilon,\nu,j}$, to prove that
$$
\int_{B'_j}e^{2\nu\psi_{\varepsilon,\nu,j}-2\nu\varphi}d\lambda=
\int_{B'_j}\Big(\sum_{k=0}^{+\infty}|f_{\nu,j,k}|^2\Big)
e^{-2\nu\varphi}d\lambda<+\infty.
$$
By the strong Noetherian property of coherent ideal sheaves
([Nar66] or [GR84]), we know that the sequence of ideal sheaves
generated by the holomorphic functions
$(f_{\nu,j,k}(z)\overline{f_{\nu,j,k}(\overline w)})_{k\le k_0}$ on
$B_j\times B_j$ is locally statio\-nary as $k_0$ increases, hence
independant of $k_0$ on $B'_j\times B'_j\compact B_j\times B_j$ for
$k_0$ large enough. As the sum of the series $\sum_k f_{\nu,j,k}(z)\overline
{f_{\nu,j,k}(\overline w)}$ is bounded by
$$\Big(\sum_k |f_{\nu,j,k}(z)|^2 \sum_k|f_{\nu,j,k}(\overline w)|^2
\Big)^{1/2}$$
and thus uniformly covergent on every compact subset of $B_j\times B_j$,
and as the space of sections of a coherent ideal sheaf is closed
under the topology of uniform convergence on compact subsets, we infer
from the Noetherian property that the holomorphic function
$\sum_{k=0}^{+\infty}f_{\nu,j,k}(z)\overline{f_{\nu,j,k}(\overline w)}$
is a section of the coherent ideal sheaf generated by
$(f_{\nu,j,k}(z)\overline{f_{\nu,j,k}(\overline w)})_{k\le k_0}$ over
$B'_j\times B'_j$, for $k_0$ large enough. Hence, by restricting
to the conjugate diagonal $w=\overline z$, we get
$$
\sum_{k=0}^{+\infty}|f_{\nu,j,k}(z)|^2\le C
\sum_{k=0}^{k_0}|f_{\nu,j,k}(z)|^2\qquad\hbox{on}~~B'_j.
$$
This implies
$$
\int_{B'_j}\Big(\sum_{k=0}^{+\infty}|f_{\nu,j,k}|^2\Big)e^{-2\varphi}d\lambda
\le C
\int_{B'_j}\Big(\sum_{k=0}^{k_0}|f_{\nu,j,k}|^2\Big)e^{-2\varphi}d\lambda
= C(k_0+1).
$$
Property (2.2.10) is proved.\medskip

\noindent{\it Step 3. Subadditivity of the approximating
sequence $\psi_{\varepsilon,\nu}$.}

\noindent
We want to compare $\psi_{\varepsilon,\nu_1+\nu_2}$ and
$\psi_{\varepsilon,\nu_1}$, $\psi_{\varepsilon,\nu_2}$
for every pair of indices~$\nu_1$, $\nu_2$, first when the functions are
associated with the same covering $X=\bigcup B_j$. Consider a function
$f\in\cH_{\nu_1+\nu_2,j}$ with
$$
\int_{B_j}|f(z)|^2e^{-2(\nu_1+\nu_2)\varphi_j(z)}d\lambda(z)\le 1,\qquad
\varphi_j(z)=\varphi(z)+c_j|z^{(j)}|^2.
$$
We may view $f$ as a function $\hat f(z,z)$ defined on the diagonal
$\Delta$ of $B_j\times B_j$. Consider the Hilbert space of holomorphic
functions $u$ on $B_j\times B_j$ such that
$$
\int_{B_j\times B_j}|u(z,w)|^2
e^{-2\nu_1\varphi_j(z)-2\nu_2\varphi_j(w)}d\lambda(z)d\lambda(w)<+\infty.
$$
By the Ohsawa-Takegoshi $L^2$ extension theorem [OT87], there exists a
function $\wt f(z,w)$ on $B_j\times B_j$ such that $\wt f(z,z)=f(z)$ and
$$
\eqalign{
\int_{B_j\times B_j}|\wt f(z,w)|^2
& e^{-2\nu_1\varphi_j(z)-2\nu_2\varphi_j(w)}d\lambda(z)d\lambda(w)\cr
&\qquad{}
\le C_7\int_{B_j}|f(z)|^2e^{-2(\nu_1+\nu_2)\varphi_j(z)}d\lambda(z)=C_7,\cr}
$$
where the constant $C_7$ only depends on the dimension $n$
(it is actually independent of the radius $r_j$ if say $0<r_j\le 1$).
As the Hilbert
space under consideration on $B_j\times B_j$ is the completed tensor
product $\cH_{\nu_1,j}\mathop{\wh\otimes}\cH_{\nu_2,j}$, we infer that
$$
\wt f(z,w)=\sum_{k_1,k_2}c_{k_1,k_2}f_{\nu_1,j,k_1}(z)
f_{\nu_2,j,k_2}(w)
$$
with $\sum_{k_1,k_2}|c_{k_1,k_2}|^2\le C_7$. By restricting to the
diagonal, we obtain
$$
|f(z)|^2=|\wt f(z,z)|^2\le\sum_{k_1,k_2}|c_{k_1,k_2}|^2
\sum_{k_1}|f_{\nu_1,j,k_1}(z)|^2\sum_{k_2}|f_{\nu_2,j,k_2}(z)|^2.
$$
>From (2.2.3) and (2.2.4), we get
$$
\psi_{\varepsilon,\nu_1+\nu_2,j}\le {\log C_7\over\nu_1+\nu_2}+
{\nu_1\over\nu_1+\nu_2}\psi_{\varepsilon,\nu_1,j}+{\nu_2\over\nu_1+\nu_2}
\psi_{\varepsilon,\nu_2,j},
$$
in particular
$$
\psi_{\varepsilon,2^\nu,j}\le \psi_{\varepsilon,2^{\nu-1},j}+{C_8\over 2^\nu},
$$
and we see that $\psi_{\varepsilon,2^\nu}+C_8 2^{-\nu}$ is a decreasing
sequence. By Step~2 and Lebesgue's monotone convergence theorem,
we infer that for every $\varepsilon,\delta>0$ and $a\le a_0\ll 0$ fixed,
the integral
$$
I_{\varepsilon,\delta,\nu}=\int_X
\Big(e^{-2\varphi}-e^{-2\max(\varphi,(1+\delta)(\psi_{2^\nu,\varepsilon}+a
))}\Big)dV_\omega
$$
converges to $0$ as $\nu$ tends to~$+\infty$ (take
$\ell={1\over\delta}+1$ and $2^\nu>\ell$ and $a_0$ such that
$\delta\sup_X\varphi+a_0\le 0$; we do not have monotonicity strictly speaking
but need only replace $a$ by $a+C_8 2^{-\nu}$ to get it, thereby slightly
enlarging the integral).
\medskip

\noindent
{\it Step 4. Selection of a suitable upper envelope.}

\noindent For the simplicity of notation, we assume here that
$\sup_X\varphi=0$ (possibly after subtracting a constant), hence we can take
$a_0=0$ in the above. We may even further assume that all our
functions $\psi_{\varepsilon,\nu}$ are nonpositive. By Step~3,
for each $\delta=\varepsilon=2^{-k}$, we can select an index
$\nu=p(k)$ such that
$$
I_{2^{-k},2^{-k},p(k)}=
\int_X\Big(e^{-2\varphi}-e^{-2\max(\varphi,(1+2^{-k})\psi_{2^{-k},2^{p(k)}})}
\Big)dV_\omega\le 2^{-k}
\leqno(2.2.11)
$$
By construction, we have an estimate $i\ddbar\psi_{2^{-k},2^{p(k)}}\ge
\gamma-C_52^{-k}\omega$, and the functions $\psi_{2^{-k},2^{p(k)}}$ are
quasi-psh with logarithmic poles. Our estimates (especially (2.2.7))
imply that $\lim_{k\to+\infty}\psi_{2^{-k},2^{p(k)}}(z)=\varphi(z)$
as soon as $2^{-p(k)}\log\big(1/\inf_j r_j(k)\big)\to 0$ (notice that
the $r_j$'s now depend on $\varepsilon=2^{-k}$). We set
$$
\varphi_\nu(z)=\sup_{k\ge\nu}(1+2^{-k})\psi_{2^{-k},2^{p(k)}}(z).
\leqno(2.2.12)
$$
By construction $(\varphi_\nu)$ is a decreasing sequence and satisfies the
estimates
$$\varphi_\nu\ge
\max\big(\varphi,(1+2^{-\nu})\psi_{2^{-\nu},2^{p(\nu)}}\big),
\qquad
i\ddbar\varphi_\nu\ge \gamma-C_52^{-\nu}\omega.$$
Inequality (2.2.11) implies that
$$
\int_X(e^{-2\varphi}-e^{-2\varphi_\nu})dV_\omega\le
\sum_{k=\nu}^{+\infty}2^{-k}=2^{1-\nu}.
$$
Finally, if $Z_\nu$ is the set of poles of $\psi_{2^{-\nu},2^{p(\nu)}}$, then
$Z_\nu\subset Z_{\nu+1}$ and $\varphi_\nu$ is continuous on~$X\ssm Z_\nu$.
The reason is that in a neighborhood of every point
\hbox{$z_0\in X\ssm Z_\nu$},
the term \hbox{$(1+2^{-k})\psi_{2^{-k},2^{p(k)}}$} contributes to
$\varphi_\nu$ only when it is larger than
$(1+2^{-\nu})\psi_{2^{-\nu},2^{p(\nu)}}$.
Hence, by the almost-monotonicity, the relevant terms of the sup
in (2.2.12) are squeezed
between $(1+2^{-\nu})\psi_{2^{-\nu},2^{p(\nu)}}$ and
$(1+2^{-k})(\psi_{2^{-\nu},2^{p(\nu)}}+C_82^{-\nu})$, and therefore
there is uniform convergence in a neighborhood of~$z_0$.
Finally, condition c) implies that
$$
\int_U|f|^2(e^{-2\varphi}-e^{-2\varphi_\nu})dV_\omega<+\infty
$$
for every germ of holomorphic function $f\in\cO(U)$ at a point~$x\in X$.
Therefore both integrals $\int_U|f|^2e^{-2\varphi}dV_\omega$ and
$\int_U|f|^2e^{-2\varphi_\nu}dV_\omega$ are simultaneously convergent
or divergent, i.e.\ $\cI(\varphi)=\cI(\varphi_\nu)$. Theorem 2.2.1 is
proved, except that $\varphi_\nu$ is possibly just continuous instead of
being smooth. This can be arranged by Richberg's regularization theorem
[Ri68], at the expense of an arbitrary small loss in the Hessian form.\qed

\claim 2.2.13. Remark|{\rm By a very slight variation of the proof, we
can strengthen condition c) and obtain that for every $t>0$
$$
\int_X(e^{-2t\varphi}-e^{-2t\varphi_\nu})dV_\omega
$$
is finite for $\nu$ large enough and converges to $0$ as $\nu\to+\infty$.
This implies that the sequence of multiplier ideals
$\cI(t\varphi_\nu)$ is a stationary decreasing
sequence, with \hbox{$\cI(t\varphi_\nu)=\cI(t\varphi)$} for $\nu$ large.}
\endclaim

\subsection{\S2.3. A Bochner type inequality}

Let $(L,h)$ be a smooth hermitian line bundle on a (non necessarily
compact) K\"ahler manifold $(Y,\omega)$.  We denote by
$|~~|=|~~|_{\omega,h}$ the pointwise hermitian norm on
$\Lambda^{p,q}T^\star_Y\otimes L$ associated with $\omega$ and $h$,
and by $\Vert~~\Vert=\Vert~~\Vert_{\omega,h}$ the global $L^2$ norm
$$
\Vert u\Vert^2 = \int_Y |u|^2 dV_\omega \qquad
\hbox{where}\quad  dV_\omega={\omega^n\over n!}
$$
We consider the $\dbar$ operator acting on $(p,q)$-forms with values in
$L$, its adjoint $\dbar_h^\star$ with respect to $h$ and the complex
Laplace-Beltrami operator $\overline{\boxsquare}_h=
\dbar\dbar_h^\star+\dbar_h^\star\dbar$.
Let $v$ be a smooth $(n-q,0)$-form with compact support
in~$Y$. Then $u=\omega^q\wedge v$ satisfies
$$
\Vert\dbar u\Vert^2+\Vert\dbar^\star_h u\Vert^2=
\Vert\dbar v\Vert^2+\int_Y\sum_{I,J}\Big(\sum_{j\in J}\lambda_j\Big)|u_{IJ}|^2
\leqno(2.3.1)
$$
where $\lambda_1\le\ldots\le\lambda_n$ are the curvature
eigenvalues of $\Theta_h(L)$ expressed in an orthonormal frame
$(\partial/\partial z_1,\ldots,\partial/\partial z_n)$ (at some fixed
point $x_0\in Y$), in such a way that
$$
\omega_{x_0}=i\sum_{1\le j\le n}dz_j\wedge d\overline z_j,\qquad
\Theta_h(L)_{x_0}=i\ddbar\varphi_{x_0}=
i\sum_{1\le j\le n} \lambda_jdz_j\wedge d\overline z_j.
$$
The proof of (2.3.1) proceeds by checking that
$$
(\dbar^\star_\varphi\,\dbar+\dbar\,\dbar^\star_\varphi)(v\wedge\omega^q)-
(\dbar^\star_\varphi\,\dbar v)\wedge\omega^q=
q\,i\ddbar\varphi\wedge\omega^{q-1}\wedge v,
\leqno(2.3.2)
$$
taking the inner product with $u=\omega^q\wedge v$ and integrating by
parts in the left hand side. In order to check (2.3.2), we use the
identity $\dbar^\star_\varphi=e^\varphi\dbar^\star(e^{-\varphi}\bu)
=\dbar^\star+\nabla^{0,1}\varphi\ort\bu$~. Let us work in
a local trivialization of $L$ such that $\varphi(x_0)=0$ and
$\nabla\varphi(x_0)=0$. At~$x_0$ we then find
$$
\eqalign{
(\dbar^\star_\varphi\,\dbar+\dbar\,\dbar^\star_\varphi)(\omega^q\wedge v)&{}-
\omega^q\wedge(\dbar^\star_\varphi\,\dbar v)={}\cr
&\big[(\dbar^\star\,\dbar+\dbar\,\dbar^\star)(\omega^q\wedge v)-
\omega^q\wedge(\dbar^\star\,\dbar v)\big]+\dbar(\nabla^{0,1}\varphi\ort(
\omega^q\wedge v)).\cr}
$$
However, the term $[\,\ldots\,]$ corresponds to the case of a trivial vector
bundle and it is well known in that case that $[\overline{\boxsquare},
\omega^q\wedge\bu]=0$, hence $[\,\ldots\,]=0$. On the other hand
$$\nabla^{0,1}\varphi\ort(\omega^q\wedge v)=
q(\nabla^{0,1}\varphi\ort\omega)\wedge \omega^{q-1}\wedge v=
-q\,i\partial\varphi\wedge \omega^{q-1}\wedge v,$$
and so
$$
(\dbar^\star_\varphi\,\dbar+\dbar\,\dbar^\star_\varphi)(\omega^q\wedge v)-
\omega^q\wedge(\dbar^\star_\varphi\,\dbar v)=q\,i\ddbar\varphi
\wedge \omega^{q-1}\wedge v.$$
Our formula is thus proved when $v$ is smooth and compactly
supported. In general, we have:

\claim 2.3.3. Proposition|Let $(Y,\omega)$ be a {\rm complete} K\"ahler
manifold and $(L,h)$ a smooth hermitian line bundle such that the
curvature possesses a uniform lower bound\break $\Theta_h(L)\ge -C\omega$.
For every measurable $(n-q,0)$-form $v$ with $L^2$ coefficients and
values in $L$ such that $u=\omega^q\wedge v$ has differentials
$\dbar u$, $\dbar^\star u$ also in $L^2$, we have
$$
\Vert\dbar u\Vert^2+\Vert\dbar^\star_h u\Vert^2=
\Vert\dbar v\Vert^2+\int_Y\sum_{I,J}\Big(\sum_{j\in J}\lambda_j\Big)|u_{IJ}|^2
$$
$($here, all differentials are computed in the sense of distributions$)$.
\endclaim

\proof. Since $(Y,\omega)$ is assumed to be complete,
there exists a sequence of smooth forms $v_\nu$ with compact support in $Y$
(obtained by truncating $v$ and taking the convolution with a regularizing
kernel) such that $v_\nu\to v$ in $L^2$ and such that
$u_\nu=\omega^q\wedge v_\nu$ satisfies $u_\nu\to u$, $\dbar u_\nu\to\dbar u$,
$\dbar^\star u_\nu\to\dbar^\star u$ in~$L^2$. By the curvature assumption,
the final integral in the right hand side of (2.3.1) must be under control
(i.e.\ the integrand becomes nonnegative if we add a term $C\Vert u\Vert^2$
on both sides, $C\gg 0$). We thus get the equality by passing to the limit
and using Lebesgue's monotone convergence theorem.\qed

\subsection{\S2.4. Proof of Theorem 2.1.1}

To fix the ideas, we first indicate the proof in the much simpler case
when $(L,h)$ is hermitian semipositive, and then treat the general
case.\medskip

\noindent
{\bf (2.4.1) Special case.} $(L,h)$ is (smooth) hermitian semipositive

Let $\{\beta\}\in H^q(X,\Omega^n_X\otimes L)$ be an arbitrary cohomology class.
By standard $L^2$ Hodge theory, $\{\beta\}$ can be represented by a smooth
harmonic $(0,q)$-form $\beta$ with values in $\Omega^n_X\otimes L$. We can
also view $\beta$ as a $(n,q)$-form with values in $L$. The pointwise
Lefschetz isomorphism produces a unique $(n-q,0)$-form $\alpha$ such
that $\beta=\omega^q\wedge\alpha$. Proposition 2.3.3 then yields
$$
\Vert\dbar\alpha\Vert^2+\int_Y\sum_{I,J}\Big(\sum_{j\in J}\lambda_j\Big)
|\alpha_{IJ}|^2=\Vert\dbar\beta\Vert^2+\Vert\dbar^\star_h \beta\Vert^2=0,
$$
and the curvature eigenvalues $\lambda_j$ are nonnegative by our assumption.
Hence $\dbar\alpha=0$ and $\{\alpha\}\in H^0(X,\Omega^{n-q}_X\otimes L)$
is mapped to $\{\beta\}$ by $\Phi^q_{\omega,h}=\omega^q\wedge\bu~$.
\medskip

\noindent
{\bf (2.4.2) General case.}

There are several difficulties. The first difficulty is that the
metric $h$ is no longer smooth and we cannot directly represent
cohomology classes by harmonic forms. We circumvent this problem by
smoothing the metric on an (analytic) Zariski open subset and by 
avoiding the remaining poles on the complement. However, some careful 
estimates have to be made in order to take the error terms into account.

Fix $\varepsilon=\varepsilon_\nu$ and let $h_\varepsilon=h_{\varepsilon_\nu}$
be an approximation of~$h$, such that $h_\varepsilon$ is smooth on
$X\ssm Z_\varepsilon$ ($Z_\varepsilon$ being an analytic subset of $X$),
$\Theta_{h_\varepsilon}(L)\ge -\varepsilon\omega$,
$h_\varepsilon\le h$ and $\cI(h_\varepsilon)=\cI(h)$. This is possible
by Theorem~2.2.1. Now, we can find a family
$$
\omega_{\varepsilon,\delta}=\omega+\delta(i\ddbar \psi_\varepsilon+\omega),
\qquad \delta>0
$$
of {\it complete K\"ahler} metrics on $X\ssm Z_\varepsilon$, where
$\psi_\varepsilon$ is a quasi-psh function on $X$ with
$\psi_\varepsilon=-\infty$ on $Z_\varepsilon$,
$\psi_\varepsilon$ on $X\ssm Z_\varepsilon$ and $i\ddbar
\psi_\varepsilon+\omega\ge 0$ (see e.g.\ [De82], Th\'eor\`eme 1.5).
By construction, $\omega_{\varepsilon,\delta}\ge\omega$ and
$\lim_{\delta\to 0}\omega_{\varepsilon,\delta}=\omega$.
We look at the $L^2$ Dolbeault complex $K^\bu_{\varepsilon,\delta}$
of $(n,\bu)$-forms on $X\ssm Z_\varepsilon$, where the $L^2$ norms are
induced by $\omega_{\varepsilon,\delta}$ on differential forms and by
$h_\varepsilon$ on elements in~$L$. Specifically
$$
K^q_{\varepsilon,\delta}=\Big\{u:X\ssm Z_\varepsilon\to\Lambda^{n,q}
T^\star_X\otimes L;\int_{X\ssm Z_\varepsilon}\kern-10pt
(|u|^2_{\Lambda^{n,q}\omega_{\varepsilon,\delta}
\otimes h_\varepsilon}+|\dbar u|^2_{\Lambda^{n,q+1}\omega_{\varepsilon,\delta}
\otimes h_\varepsilon})dV_{\omega_{\varepsilon,\delta}}<\infty\Big\}.
$$
Let $\cK^q_{\varepsilon,\delta}$ be the corresponding sheaf of germs
of locally $L^2$ sections on $X$ (the local $L^2$ condition
should hold on $X$, not only on $X\ssm Z_\varepsilon\,$!). Then,
for all $\varepsilon>0$ and $\delta\ge 0$,
$(\cK^q_{\varepsilon,\delta},\dbar)$ is a resolution of the sheaf
$\Omega^n_X\otimes L\otimes\cI(h_\varepsilon)=
\Omega^n_X\otimes L\otimes\cI(h)$. This is
because $L^2$ estimates hold locally on small Stein open sets, and the
$L^2$ condition on $X\ssm Z_\varepsilon$ forces holomorphic sections
to extend across~$Z_\varepsilon$ ([De82], Lemme 6.9).

Let $\{\beta\}\in H^q(X,\Omega^n_X\otimes L\otimes\cI(h))$ be a
cohomology class represented by a smooth form with values in
$\Omega^n_X\otimes L\otimes\cI(h)$ (one can use a \v Cech cocycle
and convert it to an element in the $C^\infty$ Dolbeault complex by
means of a partition of unity, thanks to the usual De Rham-Weil
isomorphism). Then
$$
\Vert\beta\Vert_{\varepsilon,\delta}^2\le \Vert\beta\Vert^2=
\int_X|\beta|^2_{\Lambda^{n,q}\omega\otimes h}dV_\omega<+\infty.
$$
The reason is that $|\beta|^2_{\Lambda^{n,q}\omega\otimes h}dV_\omega$
decreases as $\omega$ increases. This is just an easy calculation,
shown by comparing two metrics $\omega$, $\omega'$ which are
expressed in diagonal form in suitable coordinates; the norm
$|\beta|^2_{\Lambda^{n,q}\omega\otimes h}$ turns out to decrease
faster than the volume $dV_\omega$ increases; see e.g.\ [De82], Lemme 3.2;
a special case is $q=0$, then $|\beta|^2_{\Lambda^{n,q}\omega\otimes h}
dV_\omega=i^{n^2}\beta\wedge\overline\beta$ with the identification
$L\otimes\overline L\simeq\bC$ given by the metric $h$, hence the
integrand is even independent of $\omega$ in that case. 

By the proof of the De Rham-Weil isomorphism, the map
$\alpha\mapsto\{\alpha\}$ from the cocycle space
$Z^q(\cK^\bu_{\varepsilon,\delta})$ equipped with its $L^2$ topology,
into $H^q(X,\Omega^n_X\otimes L\otimes\cI(h))$ equipped with its
finite vector space topology, is continuous.  Also, Banach's open
mapping theorem implies that the coboundary space
$B^q(\cK^\bu_{\varepsilon,\delta})$ is closed in
$Z^q(\cK^\bu_{\varepsilon,\delta})$. This is true for all $\delta\ge
0$ (the limit case $\delta=0$ yields the strongest $L^2$ topology in
bidegree $(n,q)$). Now, $\beta$ is a $\dbar$-closed form in the Hilbert space
defined by $\omega_{\varepsilon,\delta}$ on $X\ssm Z_\varepsilon$, so there
is a $\omega_{\varepsilon,\delta}$-harmonic form
$u_{\varepsilon,\delta}$ in the same cohomology class
as $\beta$, such that
$$
\Vert u_{\varepsilon,\delta}\Vert_{\varepsilon,\delta}\le
\Vert\beta\Vert_{\varepsilon,\delta}.
$$
\claim 2.4.3. Remark|{\rm
The existence of a harmonic representative holds true only for $\delta>0$, 
because we need to have a complete K\"ahler metric on $X\ssm Z_\varepsilon$.
The trick of employing $\omega_{\varepsilon,\delta}$ instead of a fixed
metric $\omega$, however, is not needed when $Z_\varepsilon$ is (or can 
be taken to be) empty. This is the case if $(L,h)$ is such that 
$\cI(h)=\cO_X$ and $L$ is nef. Indeed, in that case, from the very 
definition of nefness, it is easy to prove that we can take the 
$\varphi_\nu$'s to be everywhere smooth in Theorem~2.2.1. However,
we will see in \S~2.5 that multiplier ideal sheaves are needed even
in case $L$ is nef, when $\cI(h)\ne\cO_X$.}
\endclaim

Let $v_{\varepsilon,\delta}$ be the unique $(n-q,0)$-form such that
$u_{\varepsilon,\delta}=v_{\varepsilon,\delta}\wedge
\omega_{\varepsilon,\delta}^q$ ($v_{\varepsilon,\delta}$ exists by the
pointwise Lefschetz isomorphism). Then
$$
\Vert v_{\varepsilon,\delta}\Vert_{\varepsilon,\delta}=
\Vert u_{\varepsilon,\delta}\Vert_{\varepsilon,\delta}\le
\Vert\beta\Vert_{\varepsilon,\delta}\le\Vert\beta\Vert.
$$
As $\sum_{j\in J}\lambda_j\ge -q\varepsilon$ by the assumption on
$\Theta_{h_\varepsilon}(L)$, the Bochner formula yields
$$
\Vert\dbar v_{\varepsilon,\delta}\Vert_{\varepsilon,\delta}^2\le
q\varepsilon\Vert u_{\varepsilon,\delta}\Vert_{\varepsilon,\delta}^2
\le q\varepsilon\Vert\beta\Vert^2.
$$
These uniform bounds imply that there are subsequences $u_{\varepsilon,
\delta_\nu}$ and $v_{\varepsilon,\delta_\nu}$ with $\delta_\nu\to 0$,
possessing weak-$L^2$ limits $u_\varepsilon=
\lim_{\nu\to+\infty}u_{\varepsilon,\delta_\nu}$
and $v_\varepsilon=\lim_{\nu\to+\infty}v_{\varepsilon,\delta_\nu}$.
The limit $u_\varepsilon=
\lim_{\nu\to+\infty}u_{\varepsilon,\delta_\nu}$ is with respect to
$L^2(\omega)=L^2(\omega_{\varepsilon,0})$. To check this, notice that
in bidegree $(n-q,0)$, the space $L^2(\omega)$ has the weakest topology
of all spaces $L^2(\omega_{\varepsilon,\delta})$; indeed, an easy calculation
as in ([De82], Lemme 3.2) yields
$$
|f|^2_{\Lambda^{n-q,0}\omega\otimes h}dV_\omega\le
|f|^2_{\Lambda^{n-q,0}\omega_{\varepsilon,\delta}\otimes h}
dV_{\omega_{\varepsilon,\delta}}\qquad
\hbox{if $f$ is of type $(n-q,0)$}.
$$
On the other hand, the limit
$v_\varepsilon=\lim_{\nu\to+\infty}v_{\varepsilon,\delta_\nu}$
takes place in all spaces $L^2(\omega_{\varepsilon,\delta})$, $\delta>0$,
since the topology gets stronger and stronger as $\delta\downarrow 0$
[$\,$possibly not in $L^2(\omega)$, though, because in bidegree $(n,q)$
the topology of $L^2(\omega)$ might be strictly stronger than that 
of all spaces $L^2(\omega_{\varepsilon,\delta})\,$].
The above estimates yield
$$
\eqalign{
&\Vert v_\varepsilon\Vert^2_{\varepsilon,0}=
\int_X|v_\varepsilon|^2_{\Lambda^{n-q,0}\omega\otimes h_\varepsilon}
dV_\omega\le\Vert\beta\Vert^2,\cr
&\Vert \dbar v_\varepsilon\Vert^2_{\varepsilon,0}\le q\varepsilon
\Vert\beta\Vert^2_{\varepsilon,0},\cr
\noalign{\vskip3pt}
&u_\varepsilon=\omega^q\wedge v_\varepsilon\equiv\beta
\qquad\hbox{in}~~H^q(X,\Omega^n_X\otimes L\otimes\cI(h_\varepsilon)).\cr}
$$
Again, by arguing in a given Hilbert space $L^2(h_{\varepsilon_0})$,
we find $L^2$ convergent subsequences $u_\varepsilon\to u$,
$v_\varepsilon\to v$ as $\varepsilon\to 0$, and in this way get
$\dbar v=0$ and
$$
\eqalignno{
&\Vert v\Vert^2\le \Vert \beta\Vert^2,\cr
&u=\omega^q\wedge v\equiv \beta
\qquad\hbox{in}~~H^q(X,\Omega^n_X\otimes L\otimes\cI(h)).\cr}
$$
Theorem 2.1.1 is proved. Notice that the equisingularity property
$\cI(h_\varepsilon)=\cI(h)$ is crucial in the above proof,
otherwise we could not infer that $u\equiv \beta$ from the fact
that $u_\varepsilon\equiv \beta$. This is true only because
all cohomology classes $\{u_\varepsilon\}$ lie in the same fixed 
cohomology group $H^q(X,\Omega^n_X\otimes L\otimes\cI(h))$, whose
topology is induced by the topology of $L^2(\omega)$ on $\dbar$-closed
forms (e.g.\ through the De Rham-Weil isomorphism).\qed

\subsection{\S2.5. A counterexample}

In view of Corollary 2.1.2, one might wonder whether the morphism
$\Phi^q_\omega$ would not still be surjective when $L$ is a nef
vector bundle. We will show that this is unfortunately not so,
even in the case of algebraic surfaces.

Let $B$ be an elliptic curve and $V$ the rank let
$2$ vector bundle over $B$ which is defined as the (unique) non split
extension
$$0 \to \cO_B \to V \to \cO_B \to 0.$$
In particular, the bundle $V$ is numerically flat, i.e.\ $c_1(V)=0$,
$c_2(V)=0$. We consider the ruled surface $X=\bP(V)$. On that surface
there is a unique section $C=\bP(\cO_B) \subset X$ with $C^2=0$ and
$$\cO_X (C)= \cO_{\bP (V)} (1)$$
is a nef line bundle. It is easy to see that
$$h^0 (X,\cO_{\bP (V)} (m))= h^0 (B,S^m V)=1$$
for all $m \in \bN$ (otherwise we would have $mC=aC+M$ where
$aC$ is the fixed part of the linear system $|mC|$ and $M\ne 0$ the
moving part, thus $M^2\ge 0$ and $C\cdot M>0$, contradiction).
We claim that
$$h^0 (X,\Omega_X^1 (kC))=2$$
for all $k \ge2.$ This follows by tensoring the exact sequence
$$
0 \to \Omega^1_{X \vert C} \to \Omega_X^1 \to \pi^* \Omega^1_{C}
\simeq \cO_{C} \to 0
$$
by $\cO_X (kC)$ and observing that
$$\Omega^1_{X \vert C}= K_X=\cO_X (-2C).$$
>From this, we get
$$
0 \to H^0(X,\cO_X((k-2)C)) \to H^0(X,\Omega_X^1\cO(kC))
\to H^0(X,\cO_X(kC))
$$
where $h^0 (X,\cO_X ((k-2)C))=h^0 (X,\cO_X (kC))=1$ for
all $k\ge 2$. Moreover, the last arrow is surjective because we
can multiply a section of $H^0(X,\cO_X(kC))$ by a nonzero section
in $H^0(X,\pi^*\Omega_B^1)$ to get a preimage. Our claim follows.
We now consider the diagram
$$
\matrix{
H^0 (X,\Omega_X^1 (2C)) & \buildo\vlra^{\wedge \omega} &
H^1 (X,K_X (2C)) \cr
\noalign{\vskip4pt}
\simeq \Big\downarrow & & \Big\downarrow \varphi \cr
H^0 (X,\Omega_X^1 (3C)) &
\buildo\vlra^{\wedge \omega}_{\psi} & H^1(X,K_X (3C)).\cr
}
$$
Since $K_X (2C) \simeq \cO_X$ and $K_X (3C) \simeq \cO_X(C),$ the
cohomology sequence of
$$
0 \to K_X (2C) \to K_X (3C) \to K_X (3C) \vert C \simeq \cO_{C}
\to 0
$$
immediately implies $\varphi=0$ (notice that $h^1 (X,K_X (2C))=h^1 (X,K_X
(3C))=1,$ since $h^1 (B,\cO_B)=h^1 (B,V)=1)$, and
$h^2(X,K_X(2C))=h^2(B,\cO_B)=0$). Therefore the diagram
implies $\psi=0$, and we get:

\claim 2.5.1.\ Proposition|$L=\cO_{\bP(V)}(3)$ is a counterample to
2.1.2 in the nef case.
\endclaim

\noindent
By Corollary 2.1.2, we infer that $\cO_X(3)$ cannot be hermitian
semi-positive and we thus again obtain -- by a quite different method --
the result of [DPS94], example 1.7.

\claim 2.5.2.\ Corollary|Let $B$ be an elliptic curve, $V$ the vector
bundle given by the unique non-split extension
$$0 \to \cO_B \to V\to \cO_B \to 0.$$
Let $X=\bP(V)$. Then $L=\cO_X(1)$ is nef but not hermitian
semi-positive $($nor does any multiple, e.g.\ the
anticanonical line bundle $-K_X=\cO_X(-2)$ is nef but not
semi-positive$)$.
\endclaim

\noindent
We now show that the above counterample is the only one that can occur
on a surface, at least when $L=\cO_X(\lambda C)$ and $C$ is an elliptic
curve with $C^2=0$ (if $C$ is a curve with $C^2>0$, $L$ is big and
therefore the conclusion is positive as well).

\claim 2.5.3.\ Proposition|Let $X$ be a smooth minimal compact K\"ahler
surface with K\"ahler form $\omega$. Let $C \subset X$ be a smooth
elliptic curve with $C^2=0.$ Then the natural map
$$\displaystyle \Phi^1_\omega: H^0 (X,\Omega_X^1 \otimes \cO_X
(\lambda C)) \to H^1 (X, K_X \otimes \cO_X (\lambda C))$$
is surjective for all $\lambda \in \bN$ with the following single
exception: $X=\bP(V)$ where $V$ is the unique
non-split extension $\displaystyle 0 \to \cO_B \to V \to \cO_B \to 0$
and $L=\cO_X(1)=\cO_X(C_0)$ where $C_0=\bP(\cO_B)\subset\bP(V)$ is
the section with $C_0^2=0$.
\endclaim

\proof. (0) First notice that by 2.1.2 the result is positive
if $\cO_X (\mu C)$ is generated by global sections for some
$\mu \in \bN$, or if $\cO_X (\mu C) \otimes G$ is
generated by global sections for some $G \in \Pic^0 (X)$.
Since $K_X \cdot C=K_C\cdot C=0$, $K_X$ cannot be ample,
hence $\kappa (X) \le 1$. If $\kappa (X)=1$, then $|m K_X|$
defines an elliptic fibration $f:X \to B$ and the equality
$K_X \cdot C=0$ implies $\dim f(C)=0$. Therefore we conclude by (0).

If $X$ is a torus or an hyperelliptic surface, we can directly apply
(0). If $X$ is a K3 surface, Riemann-Roch gives $\chi(\cO_X (C))=2,$ hence
$h^0 (\cO_X (C)) \ge 2,$ so that $\cO_X (C)$ is generated by global
sections. If $X$ is Enriques, choose a 2:1 unramified cover
$h:\tilde{X} \to X$ with $\tilde{X}$ a K3 surface. Then $h^* (
\cO_X (C))$ is generated by global sections and therefore $\cO_X
(C)$ is semi-positive, so that 2.1.2 applies again.

It remains to treat the case $\kappa (X)= -\infty.$ Since $X \not= \bP
_2,$ the surface $X$ carries a $\bP_1$-bundle structure $f:X \to B$,
and $f (C)=B$. In particular the genus $g (B) \le 1.$ We cannot have
$B=\bP_1,$ since there is no rational ruled surface $X$ carrying an
elliptic curve $C$ with $C^2=0,$ as we check immediately by [Ha77, V.2].
Hence $B$ is elliptic.
In that case [Ha77, V.2] gives immediately that $X=\bP (V)$ with a
semi-stable rank 2 vector bundle $V$ on $B$. We normalize $V$ in
such a way that $c_1 (V) \in \{0,1\}.$
\medskip

\noindent
(a) $c_1 (V)=0.$
\medskip

\noindent
Then either $V=\cO \oplus L$ with $L \in \Pic^0(B)$
or there is a non-split extension
$$0 \to \cO_B \to V \to \cO_B \to 0.\leqno(2.5.4)$$
If $V$ splits, $\cO_{\bP (V)} (1)$ is semi-positive. Since
$$\cO_X (C) \equiv \cO_{\bP (V)} (\alpha)$$
for some $\alpha \in \bN$, we conclude by 2.1.2.
In the non split case, we claim that the curve $C$ must be equal to
$C_0=\bP(\cO_B)\subset\bP(V)$. In fact, sequence $(2.5.4)$ implies that
$$H^0 (\cO_{\bP (V)} (\alpha) \otimes \pi^* (L))=0$$
for all $L \in
\Pic^0 (B), L \neq \cO_B.$ Hence
$$
\cO_X (C) \simeq \cO_X (\alpha\, C_0)$$
for some $\alpha \in \bN.$ If $\alpha \neq 1$ or $C \neq
C_0,$ then $V$ would split -- possibly after taking a finite \'etale
cover $C \to B$, which is not the case.

The latter case was already discussed and leads to the exception
mentioned in the theorem.
\medskip

\noindent
(b) $c_1 (V)=1.$

\noindent
We perform a base change $h:C \to B, \tilde{V}=h^* (V).$
Let $\tilde{X}=\bP (\tilde{V}).$ Then $\tilde{V}$ is semi-stable with
$c_1 (\tilde{V})$ even so that we are in case (a). By [At57], $C$ is
2:1 (\'etale) over $B$. Therefore $h^{-1}(C)$ consists of two
sections, hence $\tilde{V}$ splits.
Consequently we can easily reduce ourselves to the splitting case of
(a) and obtain surjectivity.\qed
\medskip

\subsection{\S2.6. A direct image theorem}
\medskip

We state here, for later use, the following useful direct image theorem.

\claim 2.6.1.\ Theorem|Let $f:X\to Y$ be a holomorphic $($smooth$)$
submersion between compact K\"ahler manifolds, and let $L$ be a
pseudo-effective line bundle.
We assume that the direct image $\cE=f_*(K_{X/Y}\otimes L)$ is locally free
and that the zero variety of $\cI(h_{\min})$ does not project onto~$Y$.
Set $\cJ=f_\star(\cI(h_{\min}))\subset\cO_Y$, $\cJ\ne 0$.
\smallskip
\item{\rm a)} If $Y$ is projective, there exists a very ample line bundle
$G$ on $Y$ such that the global sections of
$\cE^{\otimes m}\otimes G$ generate a subsheaf containing
$\cE^{\otimes m}\otimes G\otimes\cJ^m$ for every integer $m>0$.
\smallskip
\item{\rm b)} If $Y$ is projective or K\"ahler, then $\cO_{\bP(\cE)}(1)$ is
pseudo-effective on $\bP(\cE)$.
\smallskip
\item{\rm c)} If $\cE$ is of rank $1$, then $\cE$ is a pseudo-effective
line bundle on~$Y$.
\vskip0pt
\endclaim

\proof. The proof closely follows ideas already described by Viehweg,
[Kol86], [DPS94] and [Mou97], so we will be rather quick on details.
\medskip

\noindent
a) One can take for instance $G=K_Y+(n+1)A$ where $A$ is very ample on $Y$
and $n=\dim Y$. Then
$$
K_{X/Y}\otimes L\otimes f^\star G=K_X\otimes L\otimes f^\star A
$$
and $L\otimes f^\star A$ can be equipped with the tensor product of
the metric $h_{\min}$ of $L$ by a singular metric on $A$ which is
smooth of positive curvature outside a point $y\in Y$, with a single
isolated pole of Lelong number${}>1-\varepsilon$ at $y$ ([De90], \S$\,$6).
H\"ormander's standard $L^2$ estimates show that sections of
$K_{X/Y}\otimes L\otimes f^\star G\otimes\cI(h_{\min})$ on $X_y=f^{-1}(y)$
can be extended to global sections on $X$; actually,
given such a section $h$ defined on $f^{-1}(V)$, where $V$ is a
neighborhood of~$y$, we solve the $\dbar$ equation
$\dbar u=\dbar(\theta(f)h)=h\,f^\star\dbar\theta$ where
$\theta$ is a cut-off function with support in $V$, equal to $1$
near $y$. By construction, the curvature current of
$L\otimes f^\star G$ satisfies
$\Theta\ge f^\star\omega_Y$ for some K\"ahler form $\omega_Y$ on $Y$.
The curvature need not be positive on $X$, but this is nevertheless
sufficient to solve the $\dbar$-equation in virtue of ([De82], Th\'eor\`eme
4.1), since the norm $|f^\star\dbar\theta|_{\Theta}$ is bounded
(in the notation of [De82]). Moreover, the
Lelong number of the induced metric on $L\otimes f^\star G$ along
the fiber $X_y$ will be $(n+1)(1-\varepsilon)$, thus in the range
$]n,n+1[$, so that the resulting Nadel
multiplier ideal sheaf $\cI'$ of that metric satisfies $\cI'\subset
\cI(h_{\min})\cap\cI_{X_y}$ by [Sk72]. This implies that the
solution $u$ vanishes along $X_y$ and that
$$
\wt h=\theta h-u\in H^0(X,K_{X/Y}\otimes L\otimes f^\star G\otimes
\cI(h_{\min}))
$$
coincides with $h$ in restriction to $X_y$. In other words, the direct
image sheaf
$$
f_*\big(K_{X/Y}\otimes L\otimes f^\star G\otimes\cI(h_{\min})\big)
$$
is generated by global sections. However, this sheaf is obviously contained
in $\cE\otimes G$ and its global sections contain those of
$\cE\otimes G\otimes\cJ$. The assertion for $\cE^m$ follows by
the usual fiber product trick, where $X\to Y$ is replaced by
$X_m=X\times_Y\ldots\times_Y X\to Y$ (recall that $f:X\to Y$ is
supposed to be smooth). Then $K_{X_m/Y}\otimes
(L\stimes_Y\ldots\stimes_Y L)$ has direct image $\cE^{\otimes m}$ on $Y$.
\medskip

\noindent
b) is an straightforward consequence of a), at least in the projective
situation, since $G$ gets multiplied by $1/m$ as $m$ goes to~$+\infty$.
The K\"ahler case (which we will not need anyway) can be dealt with as
in Mourougane [Mou97], by using metrics and local sections over a
fixed Stein covering of $Y$.

\noindent
c) special case of b).\qed

\subsection{\S2.7. Applications}

Our applications mostly concern compact K\"ahler manifolds such that
either the canonical or anticanonical line bundle is pseudo-effective.
The first one has been observed independently by M.\ Paun [Pa98].

\claim 2.7.1.\ Proposition| Let $X$ be a compact K\"ahler manifold with
$-K_X$ pseudo-effective. Assume that $-K_X$ has a $($singular$)$
hermitian metric $h$ with semi-positive curvature such that
$\cI(h)=\cO_X$ $($i.e.\ the singularities of the weights $\varphi$
are mild enough to warrant that $e^{-\varphi}$ is locally integrable$)$.
Then
\smallskip
\item{\rm a)} The natural pairing
$$H^0(X,T_X)\times H^0(X,\Omega^1_X)\to \bC$$
is non degenerate on the $H^0(X,\Omega^1_X)$ side,
and the non zero holomorphic $1$-forms do not vanish at all.
\smallskip
\item{\rm b)} The Albanese map $\alpha:X\to\Alb(X)$ is a submersion, and there
is a group of automorphisms of $X$ which lies above the translations of the
Albanese torus.
\vskip0pt
\endclaim

\proof. a) The hard Lefschetz theorem applied with $L=K_X^{-1}$ and $q=1$
implies that there is a surjective map
$$
H^0(X,\Omega^{n-1}_X\otimes K_X^{-1})\buildo\vlra^
{\omega\wedge\bu} H^1(X,K_X\otimes K_X^{-1})=H^1(X,\cO_X).
$$
However, $H^0(X,\Omega^{n-1}_X\otimes K_X^{-1})\simeq H^0(X,T_X)$ and
the arrow $\omega\wedge\bu$ can then be seen as the contraction
mapping $\xi\mapsto \xi\ort\omega$ of the K\"ahler form by a
holomorphic vector field~$\xi$. Since the group $H^1(X,\cO_X)$ is conjugate
to $H^0(X,\Omega^1_X)$ by Hodge symmetry, the (non degenerate) $L^2$ pairing
between the $(0,1)$-class $\{\xi\ort \omega\}\in H^1(X,\cO_X)$ and a form
$\eta\in H^0(X,\Omega^1_X)$ is given by
$$
\int_X\langle \xi\ort \omega, \overline\eta\rangle_\omega\omega^n=
\int_X\langle \xi,\eta\rangle_{T_X\times\Omega^1_X}\omega^n=
C\langle \xi,\eta\rangle_{T_X\times\Omega^1_X},\qquad C>0.
$$
($\langle\xi,\eta\rangle_{T_X\times\Omega^1_X}$ is a holomorphic
function, hence constant).
Because of surjectivity onto $H^1(X,\cO_X)$, there exists for every
nonzero holomorphic $1$-form $\eta$ a vector field $\xi$ such that
$\langle\xi,\eta\rangle_{T_X\times\Omega^1_X}\ne 0$. This implies that
$\eta$ does not vanish and a) is proved.

\noindent b) Let $u_1,\ldots,u_q$ be a basis of $H^0(X,\Omega^1_X)$.
The Albanese map is given by
$$
\alpha:X\to\Alb(X),\qquad x\mapsto\alpha(x)=\Big(\int_{x_0}^x u_j
\Big)_{1\le j\le q}~~\hbox{modulo periods}.
$$
Hence $d\alpha\simeq(du_1,\ldots,du_q)$, and we know by a) that
$du_1(x),\ldots,du_q(x)\in T_{X,x}^\star$ are linearly independent
at every point. This means that $\alpha$ has maximal rang $q$ at every
point, i.e.\ is a submersion. The existence of vector fields shown in
a) easily imply the assertion on automorphisms.\qed

\noindent We next consider the case when $K_X$ is pseudo-effective.

\claim 2.7.2.\ Abundance conjecture|If $X$ is a compact
K\"ahler manifold with $K_X$ pseudo-effective, then the Kodaira
dimension $\kappa(X)$ is non-negative, i.e.\ there exist non trivial
sections of $H^0(X,mK_X)$ for some $m>0$.
\endclaim

\noindent The abundance conjecture is presently known only in the
projective case, and even then, only for $\dim X\le 3$. What we can prove
from our hard Lefschetz theorem is the following partial result
in the K\"ahler case.

\claim 2.7.3.\ Theorem|Let $X$ be a compact K\"ahler manifold with
$K_X$ pseudo-effective. Assume that $K_X$ has a singular hermitian
metric of non-negative curvature, possessing ``algebraic singularities''
in the following sense\/$:$ there exists a modification
$\mu:\tilde X\to X$ such that the pullbacks of the local
plurisubharmonic weights $\varphi$ take the form
$$
\varphi\circ\mu=\sum_j\lambda_j\log|g_j|+O(1)
$$
where $O(1)$ is a bounded term, $D_j=\{g_j=0\}$ is a family of
normal crossing divisors in $\tilde X$ and $\lambda_j$ are positive rational
numbers. Then $X$ satisfies at least one of the following two properties:
\smallskip
\item{\rm a)} $\chi(X,\cO_X)=\chi(X,K_X)=0$ and there exists a nonzero
holomorphic $p$-form in $H^0(X,\Omega^p_X)$ for some odd integer $p$.
\smallskip
\item{\rm b)} There exists $q=0,1,\ldots,n$ and infinitely many positive
integers $m$ such that
$$H^0(X,\Omega^q_X\otimes\cO(mK_X))\ne 0.$$
\vskip0pt
\endclaim

\proof. Observe that sections of $H^0(X,\Omega^q_X\otimes\cO(mK_X))$ are
bimeromorphic invariants, hence we can assume that $\tilde X=X$.
Suppose that b) fails, i.e.\ that there is $m_0>0$ such that
$H^0(X,\Omega^q_X\otimes\cO(mK_X))=0$ for all $q$ and $m\ge m_0$. Then a
fortiori
$$H^0(X,\Omega^q_X\otimes\cO(mK_X)\otimes\cI(h^{\otimes m}))=0$$
and the hard Lefschetz theorem implies that
$$H^q(X,K_X\otimes\cO(mK_X)\otimes\cI(h^{\otimes m}))=0\qquad
\hbox{for all $q$ and $m\ge m_0$}.
$$
Therefore $\chi(X,\cO((m+1)K_X)\otimes\cI(h^{\otimes m}))=0$ for
$m\ge m_0$. However, the assumption on the singularities of $h$ shows that
$\cI(h^{\otimes m})=\cO_X(-[m\lambda_j]D_j)$ where $[~~]$ denotes the
integral part. By Riemann-Roch, we have
$$
\eqalign{
f(m)&=\chi(X,\cO((m+1)K_X)\otimes\cI(h^{\otimes m}))\cr
&=\int_X \exp\Big((m+1)c_1(K_X)-\sum[m\lambda_j]c_1(D_j)\Big)\Todd(X).\cr}
$$
By expanding the integral, we find a (constant) integer $N$ such that
$$
N\,f(m)=P(m,[m\lambda_1],\ldots,[m\lambda_r])
$$
for some polynomial $P\in \bZ[t_0,t_1,\ldots,t_r]$ of degree${}\le n$.
Take $m$ to be a large multiple $kd$ of a common denominator of the
$\lambda_j$'s. Then $f(kd)$ is a polynomial in $k$ and vanishes
for $k$ large, thus $f(0)=\chi(K_X)=0$. Therefore
$\chi(\cO_X)=(-1)^n\chi(K_X)=0$, and as $h^0(X,\cO_X)=1$ we conclude
that one of the odd degree groups $H^p(X,\cO_X)$ must be nonzero.
By Hodge symmetry, we get $H^0(X,\Omega^p_X)\ne 0\,$; property
a) is proved.
\medskip


It seems likely that in Theorem 2.7.2 (at least in case b), the Kodaira
dimension $\kappa(X)$ should be non-negative. For our purposes it
suffices to have

\claim 2.7.4.\ Proposition| Let $X$ be a compact manifold, $E$ a vector
bundle and $L$ a line bundle on $X.$ Suppose that $H^0(X,E \otimes
L^m) \ne 0$ for infinitely many $m.$ Then $a(X) \geq 1$ or $\kappa (L)
\geq 0.$ \endclaim

\proof. By our assumption we have inclusions $L^{-m}\to E$. Considering
the smallest subsheaf in $E$ containing all the images
and taking determinants, we obtain a line bundle $F$ and infinitely
many inclusions $L^{-m}\to F$. So $H^0(X,F \otimes L^m)
\ne 0$ for infinitely many~$m$. If $X$ carries infinitely many
irreducible hypersurfaces, then \hbox{$a(X) > 0$}
([Kra75], see also [FiFo79]).  So suppose that $X$ has only finitely
many irreducible hypersufaces $Y_i.$ Then
consider the cone $K$ generated by the $Y_i$, say in $\Pic(X)\otimes\bR.$
Then $F \otimes L^m \in K,$ hence $L \in K,$ which implies
$\kappa (L) \geq 0.$ \qed

\claim 2.7.5.\ Theorem| Let $X$ be a smooth compact K\"ahler threefold
with $K_X$ pseudo-effective. Assume that $K_X$ has a singular metric
as in Theorem 2.7.3, e.g.\ that $K_X$ is hermitian semi-positive. Then
$\kappa (X) \geq 0.$ \endclaim

\proof. Since $K_X$ is pseudo-effective, $X$ cannot be uniruled.
Therefore the main result in [CP00] implies that $\kappa (X) \geq 0$
unless possibly if $X$ is simple (and $\kappa (X) = - \infty$), which
means that there is no positive-dimensional subvariety through the
very general point of $X.$ From 2.7.3 we obtain that
$$
\chi(X,\cO_X) = 0$$
or that
$$
H^0(X,\Omega^q_X \otimes \cO_X(mK_X)) \ne 0 $$
for infinitely many
$m.$ In the first case we have a $3$-form (so that $\kappa (X) \geq 0$
or a $1-$form, so that we have a non-trivial Albanese. $X$ being
simple, the Albanese map must be generically finite over the Albanese
torus. But then $\kappa (X) \geq 0.$ In the second case we conclude by
Proposition 2.7.4. \qed

\section{\S3. Pseudo-effective versus almost nef line bundles}

In this section we study pseudo-effective line bundles and their
``numerical'' counterparts, which we call almost nef line bundles.

\claim 3.1.\ Definition|Let $X$ be a projective manifold, $L$ a line
bundle on $X$. The bundle $L$ is almost nef if and only if there is a
family $A_i \subset X$, $i\in\bN$, of proper algebraic subvarieties
such that $L\cdot C \ge 0$ for all irreducible curves $C \not\subset
{\bigcup\limits_{i}} A_i$. The Zariski closure of the union of all
curves $C$ with $L\cdot C<0$ will be called the {\it non-nef locus}
of $L$.
\endclaim

\claim 3.2.\ Remark|{\rm We say that $(C_t)_{t \in T}$ is a covering
family of curves on $X$ if $T$ is compact, $X= {\bigcup\limits_{t
\in T}} C_t,$ and if $C_t$ is irreducible for general $t \in T.$}
\endclaim

With this notation, $L$ is almost nef if and only if $L\cdot C_t \ge 0$ for
all covering families $(C_t)$ of curves. Indeed, one direction is
clear, the other is an obvious Hilbert scheme argument.

\claim 3.3.\ Proposition|Let $X$ be a projective manifold. Assume
that $L$ is pseudo-effective on $X$. Then $L$ is almost nef.
\endclaim

\proof. This follows from [DPS96a, 4.3] but for the convenience of the
reader we give here a proof using directly the definition.  Choose $A$
ample such that $mL+A$ is effective for $m \ge m_0$ sufficiently
divisible. Denote this set of $m$'s by $M$. Then we can write
$$mL+A=E_m$$
with an effective divisor $E_m \subset X.$ Hence $E_m\cdot C
\ge 0$ for all $C \not\subset E_m$ and therefore $(mL+A)\cdot C \ge 0$ for $C
\not\subset E_m.$ Consequently we have
$$(mL+A)\cdot C \ge 0$$
for all $m \in M$ and all $C \not\subset {\bigcup\limits_{m \in M}} E_m.$
Hence $L\cdot C \ge 0$ for all those $C.$ \qed

\claim 3.4.\ Problem|Let $X$ be a projective manifold and $L$ an almost
nef divisor on $X$.  Is $L$ pseudo-effective?
\endclaim

\claim 3.5.\ Comments|{\rm This is in general a very hard problem (maybe even
the answer is negative).  We here point out some circumstances when
(3.4) has a positive answer.}
\endclaim

\noindent
a) $\dim X=2$ and $L$ is arbitrary.

This is already observed in [DPS96a, 4.5]. The reason is simply that
the cone of effective divisors is the cone effective curves, hence,
by dualizing, the ample cone is the dual cone to the cone of effective
divisors.  So $L \in \bar{K}_{\eff} (X)$ if and only if $L\cdot C \ge
0$ for all $C \subset X$ with $C^2 \ge 0.$
\medskip

\noindent
b) Consider now the case $L=K_X.$

Notice first that $K_X$ almost nef just says that $X$ is not uniruled.
In fact, if $X$ is uniruled, we have a covering family $(C_t)$ of
rational curves $(C_t)$ with $K_X\cdot C_t < 0,$ and conversely, if there
is a covering family $(C_t)$ with $K_X\cdot C_t < 0,$ then $X$ is uniruled
by [MM86].

Now suppose that $K_X$ is almost nef and $\dim X=3.$ Then $X$ has a minimal
model $X'$, i.e.\ $K_{X'}$ is nef.  By abundance, $\kappa (X)= \kappa
(X') \ge 0,$ in particular $K_X$ is pseudo-effective. Notice that this
is more than we asked for because a priori $K_X$ could be
pseudo-effective and $\kappa (X)=- \infty.$

In order to prove (3.4) in the case $L=K_X$ and $\dim X=n,$ we will
``only'' need the existence and finiteness of flips but we can avoid
the use of the abundance conjecture. Since $X$ is not uniruled, $X$
has a birational model $X'$ with $K_{X'}$ nef, therefore $K_{X'}$ is
pseudo-effective. In order to see that $K_X$ itself is pseudo-effective,
take a divisor $A'$ on $X'$ such that $H^0(m
K_{X'}+ A') \not= 0$ for $m \gg 0$ and sufficiently divisible.  Then
we only need to consider the two following situations

\noindent
(3.5.1) $\lambda: X \to X'$ is a divisorial contraction

\noindent
(3.5.2) $\lambda : X \rightharpoonup X'$ is of flipping type.

\noindent
In case (3.5.1) $K_X= \lambda^* (K_{X'}) + \mu E,$ where $E$ is the exceptional
divisor and $\mu >0,$ therefore clearly $K_X$ is pseudo-effective.
In case (3.5.2) let $A$ be the strict transform of $A'$ in $X$ and,
$\lambda$ being an isomorphism in codimension $1$, we have by the
Riemann extension theorem
$$H^0 (X, mK_X+A)= H^0(X',mK_{X'} + A'),$$
thus $H^0 (X, mK_X+A) \ne 0$ for $m\gg 0$ sufficiently divisible.
There is a slight difficulty: $A$ is a priori only a Weil divisor, but
since $X$ is $\bQ$-factorial, we find $\lambda$ such that $\lambda A$
is Cartier and moreover
$$H^0 (X, \lambda m K_X + \lambda A) \not= 0.$$
In total the flip conjectures imply (3.4) for $L=K_X$
(in any dimension).\qed

\section{\S4. Varieties with pseudo-effective anticanonical bundles}

In this section we study compact K\"ahler manifolds and projective varieties
with pseudo-effective and nef anticanonical bundles. We shall begin with
the nef case, in which already quite some results have been obtained.
In fact, concerning the structure of compact K\"ahler manifolds $X$ with $-K_X$
nef, we have the following

\claim 4.1.\ Conjecture {\rm ([DPS93, 96b])}|Let $X$ be a compact K\"ahler
manifold with $-K_X$ nef. Then
\smallskip
\item{\rm a)} If $\varphi:X \to Y$ is a surjective map to the normal complex
  space $Y$, then $\kappa (Y) \le 0$ where $\kappa(Y)= \kappa (\hat Y)$,
  $\hat Y$ a desingularisation).
\smallskip
\item{\rm b)} The Albanese map $\alpha:X \to \Alb (X)$ is a surjective
  submersion.
\smallskip
\item{\rm c)} $\pi_1 (X)$ is almost abelian, i.e.\ abelian up to a subgroup
  of finite index.
\endclaim

\claim 4.2.\ Remark|{\rm The status of the conjecture is as follows

\smallskip\noindent
(4.2.1) a), b) and c) hold if $-K_X$ is hermitian
semi-positive ([DPS93, DPS96b])

\smallskip\noindent
(4.2.2) a) holds and therefore surjectivity of $\alpha$, if $X$
  is projective [Zh96]

\smallskip\noindent
(4.2.3) a), b) and c) hold if $X$ is a projective 3-fold ([PS97]).

\smallskip\noindent
(4.2.4) a compact K\"ahler $n$-fold, $n \leq 4$, does not admit a
surjective map to a normal projective variety of general type
and therefore a) holds for $n$-folds, $n \leq 4.$ [CPZ98].

\smallskip\noindent
(4.2.5) $X$ does not admit a map to a curve $C$ of genus $g(C) \ge 2$
[DPS93].}
\endclaim

Our main aim is here to show that the Albanese map of a compact K\"ahler
threefold is a surjective submersion also in the non-algebraic case.

\claim 4.3.\ Theorem|Let $X$ be a compact K\"ahler 3-fold with $-K_X$
nef. Then the Albanese map $\alpha:X \to \Alb (X)$ is surjective and
has connected fibers.
\endclaim

\proof. We may assume $X$ non-algebraic (4.2) and $\kappa(X)=
-{\infty}$ [Bea83]. Let $Y= \alpha(X)$ and assume $Y \not= \Alb(X)$.
By (4.2.4) we have $\dim Y=2$. Let $\hat Y \to Y$ be a
desingularisation. Then $\kappa (\hat Y) > 0$ and $q (\hat Y) \ge 3$
by [Ue75]. If necessary, substitute $\alpha$ by its Stein
factorisation. The general fiber $F$ of $\alpha:X \to Y$ must be a
smooth rational curve; otherwise $\kappa (X) \ge 0$ by $C_{3,1}$ ([Ue87]).
Next observe that the algebraic dimension $a(Y) \le 1$. In fact
otherwise $Y$ would be Moishezon and induces 2 independent meromorphic
functions on $X$. But then $X$ is clearly algebraically connected, i.e.\ any
two points can be joined by a finite union of irreducible compact
curves. Hence $X$ is Moishezon (and therefore projective) by [Cam81].
So $a(Y) \le 1$; in particular $\kappa
(\hat Y)=1=a(Y)$. Let $\hat f : \hat Y \to C$ be the algebraic reduction
(${}={}$Iitaka fibration). Then clearly $\kappa (C) = \kappa (\hat Y) = 1,$
and $C$ is of general type. Therefore the induced meromorphic map
$X \to Y$ is actually holomorphic, contradicting (4.2.5).
It remains to prove that $\alpha$ has connected fibers. If dim $A=1$,
this is a general fact [Ue75]. If dim $A=2$, we argue as follows. If
$a(A)=2$, then, using again $C_{3,1}$, $X$ is projective and we refer
to (4.2.3). If $a(A)=1$, consider the algebraic reduction $f: A \to
B$ to the elliptic curve. If $\alpha$ is not connected, then so does
$\beta: X \to B$. Let $\gamma :X \to \tilde {B}$ be the Stein
factorisation of $\beta$. By (4.2.4), $\tilde{B}$ must be an elliptic
curve. This contradicts clearly the universal property of the Albanese
torus.  If finally $a(A)=0$, then consider the Stein factorisation
$g:X \to S$ of $\alpha$. Since $A$ contains no curves, the map $h:S
\to A$ is unramified, hence $S$ is a torus itself.  Therefore $h=\Id$. \qed

In order to investigate further the structure of compact K\"ahler 3-folds with
$-K_X$ nef, we quote the

\claim 4.4. Proposition| Let $Z$ be a compact K\"ahler threefold
and $f: Z \la C$ be a surjective map with connected fibers to a smooth
curve of genus $g \geq 1.$ Assume that the general fiber $F$ has
$\kappa (F) = - \infty$ and $q(F) = 1.$ If $-K_{Z \vert C}$ is nef,
then the only singular fibers of $f$ are multiples of smooth surfaces.
In particular, $f$ is smooth after a finite base change $\tilde C \la
C.$ If $C$ is elliptic, the original $f$ is already smooth.
Furthermore there is a smooth minimal surface $Y$, a $\bP_1$-bundle
structure $g: Z \la Y$ and an elliptic fibration $h: Y \la C$ with at
most multiple fibers as singular fibers such that $f= h \circ g.$
\endclaim

The proof can be found in [CPZ98]. For the convenience of the reader, we
sketch the idea of the proof. Since $Z$ is K\"ahler and $f$ is locally
Moishezon, $f$ is projective [CP00].  Since $K_Z$ is not $f-$nef, a
theorem of Kawamata allows us to construct locally over the base $C$ a
relative contraction, and it turns out that the dimension of the
images is always 2 unless $f$ is birational.  Then one proves that all
these local relative contractions glue to a global relative
contraction $g: Z \la Y.$ In the birational case one has to repeat
this construction and one finally ends up with some fibration
analogous to $g$ which has to be studied with the methods of [PS97].

\claim 4.5.\ Corollary| Let $X$ be a compact K\"ahler threefold with $-K_X$
nef. If $q(X)=1$, $X$ is projective, unless $K_X \equiv 0$.
\endclaim

\proof. Let $\alpha: X \to A$ be the Albanese to the elliptic curve
$A$ and suppose $X$ not projective and $K_X \not \equiv 0.$
By $C_{3,1},$ we have $\kappa (F) = - \infty$ for the general fiber $F$ of
$\alpha.$ Since $-K_F$ is nef, the irregularity $q(F) \leq 1.$ If $q(F) = 0,$
then $X$ is algebraic; see [CP00] for that and further references. Actually
we can also conclude as follows. Let $\omega$ be a non-zero 2-form. Consider
the exact sequence
$$ 0 \la N^*_F \otimes \Omega^1_F \la \Omega^2_X \vert F \la \Omega^2_F \la 0;$$
then $\omega \vert F$ induces a 1-form on $F.$
So by (4.4), $\alpha$ is a submersion
and there is a smooth surface $Y$ and submersions $g: X
\to Y$ and $h: Y\to A$ such that $f = h \circ g. $ Moreover $g$ is
$\bP_1$-bundle and $h$ is a smooth elliptic fibration. Since $A$ is an
elliptic curve, $h$ is locally trivial and therefore $\kappa (Y) = 0;$
actually $Y$ is a torus or hyperelliptic. However $Y$ cannot be
projective, otherwise $X$ is projective, so $Y$ is a torus. This
contradicts $q(X) = 1.$ \qed

Continuing the study of non-algebraic compact K\"ahler threefolds $X$
with $-K_X$ nef and $K_X \not\equiv 0$, we therefore either have
$q(X)= 0$ or $q(X)=2$.

\claim 4.6.\ Theorem|Let $X$ be a compact K\"ahler 3-fold with $-K_X$
nef, $K_X \not\equiv 0$ and $q(X)=2$. Let $\alpha:X \to A$ be the
Albanese map and assume $a(A) >0$.  Then $X$ is a $\bP_1$-bundle over
$A$.
\endclaim

\proof. By (4.3), it is clear that dim $A=2$. Note that by $(C_{3,2})$ the
general
fiber of $\alpha$ must be $\bP_1$. Therefore $a(A)=2$ implies $a(X)=3$
and $X$ is projective. We conclude by (4.2.3). So we may suppose $a(A)=1$;
let $\pi: A \la B$ be the algebraic reduction, an elliptic bundle.
Now we can apply (4.4) to the composite map $X \la B$ to obtain our claim. \qed

\noindent
We now investigate the case $a(A)=0.$

\claim 4.7.\ Theorem|Let $X$ be a compact K\"ahler 3-fold with $-K_X$
nef, $K_X \not \equiv 0$ and $q(X)= 2.$ Let $ \alpha : X \to A$ be the
Albanese and assume
$a(A)=0$. Then $\alpha$ is a projectively flat ${\bP}_1$-bundle.
\endclaim

We already know in the situation of 4.7 that $\alpha: X \la A$ is
surjective with connected fibers and the general fiber is $\bP_1.$
Let us first see that $\alpha$ is projectively flat once we know
that $\alpha $ is smooth. In fact, the exact sequence
$$ 0 \la T_{X/A} \la T_X \la T_A $$
and the observation $T_{X/A} = -K_X$ show that $T_X$ is nef. Then
the claim follows from [CP91,8.2] (the non-algebraic case is just the
same).
\smallskip \noindent Now we show that $\alpha$ is a ${\bP}_1-$bundle.
Since $A$ does not contain subvarieties of positive dimension, $\alpha
$ can have at most a finite number of singular fibers; therefore there
is a finite set $E \subset A$ such that $\alpha$ is smooth over $A_0 =
A \setminus E.$ We let $X_0 = \alpha^{-1}(A_0)$ and we must prove $E =
\emptyset.$

We prepare the proof of (4.7) by the following

\claim 4.8.\ Lemma|In the situation of (4.6) we have
$$H^0(X,-K_X \otimes \alpha^*(L))=0$$
for all $L \in \Pic(A),$ unless $X$ is a $\bP_1$-bundle.
\endclaim

\proof. Suppose $H^0(X,-K_X \otimes \alpha^*(L)) \neq 0$ and
take $D \in |-K_X \otimes \alpha^*(L)|$.  Let $D_0$ be a component of
$D$ with $\alpha(D_0)=A$. Then $\kappa(D_0) \ge 0.$
Let $\lambda_0$ be the multiplicity of $D_0$ in $D$.  Then
$$K_ { \lambda_{0} D_0}=K_D \vert \lambda_{0} \, D_0 - E$$
with $E$ effective. By adjunction
$$\lambda_0 K_{{D}_{0}}= K_{\lambda_{0} D_{0}} \vert D_0 +
(\lambda_0-1)K_X \vert D_0.$$
Taking into account $K_D \vert D_0 = \alpha^*(L) \vert D_0$ (again by
adjunction),
we obtain
$$K_{{D}_{0}} \equiv -B_1-B_2 +  \alpha^*(L')$$
with $B_1$ effective and $B_2$ nef,
where $L'$ is a rational multiple of $L.$
Let ${\tau} : \hat{D}_0 \to D_0$ be a desingularization of $D_0$;
namely a minimal desingularization after normalization. Then still
$$K_{\hat{D}_0} \equiv -{\hat{B}_1} - {\hat{B}_2 + \hat \alpha^*(L')}$$
with $\hat{B}_1$
effective, $\hat{B}_2$ nef and $\hat \alpha$ is the induced map.
But $\kappa(\hat{D}_0) \ge 0$. So $\hat
\alpha^*(L')$ is pseudo-effective; on the other hand $L'$ is of signature
$(1,1)$ if not numerically trivial ([LB92,p.318,Ex.8).  This implies
$K_{\hat{D}_0} \equiv L  \equiv 0$, and $\hat B_1 = \hat B_2 = 0,$ hence
$\hat{D}_0$ is a torus and
${\hat{D}_0} \to A$ is \'etale.
>From $\hat{B}_1 = \hat{B}_2 = 0$ we conclude that $B_1=B_2=0$ and $D_0$ is
normal with at most rational double points. Since $\hat{D}_0$ is a
torus, it has no rational curves, hence $D_0$ is smooth.  Write
$$-K_X \equiv \sum_{i \ge 0} \lambda_i D_i + \sum \mu_j R_j$$
with $\alpha(D_i)=A$, and $R_j \cdot F=0$, $F$ the general fiber of $\alpha$.
Since $a(A)=0$, we have
$$\dim\alpha (R_j)=0$$
for all $j$.  Moreover we know already $D_i\cap R_j = \emptyset$ for all
$i,j$.  Hence $\sum \mu_j R_j$ must be nef.
In fact, first observe that
$$ -K_X \vert \sum \mu_j R_j = \sum \mu_j R_j \vert \sum \mu_j R_j $$
is nef (i.e.\ nef on the reduction); then apply
[Pet98, 4.9]). But $\sum \mu_j R_j$
cannot be nef, unless all $\mu_{j}=0.$ Hence $-K_X \equiv \sum \lambda_j D_j.$

Since $K_X\cdot F=-2,$ we observe by the way $-K_X \equiv D_0+D_1$
or $-K_X \equiv 2 D_0.$ We claim that $\alpha$ is
a submersion. Suppose that $G \subset \alpha^{-1} (a)$ is a 2-dimensional fiber
component. Then we deduce, $\alpha
\vert D_i$ being \'etale, that $D_i \cap G$ is at most finite, hence
empty. From $$-K_X \equiv \sum \lambda_i D_i$$  we deduce $K_X\cdot G=0.$
Hence there must be a 1-dimensional irreducible
component $C_0 \subset \alpha^{-1} (a)$ with $D_0\cdot C_0 \not= 0.$
In particular $K_X\cdot C_0 < 0.$
In case that $\dim \alpha^{-1}(a) = 1,$ this is clear anyway, so that $C_0$
always exists. Consequently [Kol96]
$C_0$ moves in an at least 1-dimensional family, say $(C_t)_{t\in T}.$
Obviously the general $C_t$ is a fiber of $\alpha$ and
$\dim \alpha (C_t)=0$ for all $t$. So $C_0$ actually moves in a
2-dimensional family and $K_X\cdot C_0=-2.$
Now consider the graph of
$(C_t)$ and it follows immediately that $\alpha$ is a $\bP_1$-bundle. \qed
\medskip

\noindent
{\it Proof of 4.7.} The proof of 4.7 will now be completed by proving

\bigskip \noindent
{\bf 4.8.1  Claim} There exists a line $L$ on $A$ such that
$$ H^0(A,\alpha_*(-K_X) \otimes L) \ne 0.$$

\bigskip \noindent There is a slight technical difficulty: a priori
$\alpha_*(-K_X)$ need not be locally free. Therefore we will consider
its dual $W$, compare the cohomology of $W$ and $\alpha_*(-K_X)$ via
several direct image calculations and then prove equality of both
sheaves. Then it will be easy to conlcude.
To begin with, consider the exact sequence
$$0 \to \alpha_* (-K_X) \to W \to Q \to 0$$
where $W=\alpha_*(-K_X)^{**}$ (a rank 3 vector bundle) and $Q$ is just the
cokernel. Hence
$$\chi(W)= \chi( \alpha_*(-K_X)) + \chi(Q). \leqno{(4.8.1)}$$
We claim that
$$c_1(W)=0. \leqno{(4.8.2)}$$
This is seen as follows.
The exponential sequence plus the vanishing $H^q(A_0,{\cal O}) = 0$
for $q \geq 0$ (recall that $A_0$ is not compact since $E \ne \emptyset $)
yields
$$H^2(A_0, \cO_{A_0}^* ) \simeq H^3(A_0,\bZ),$$
>From that we see easily - using e.g. Mayer-Vietoris
that $H^2(\cO_{A_0}^*)$ is torsion free. Now the obstruction for a
projective bundle to come from a vector bundle is a torsion element
in $H^2(A_0,\cO^*).$ Therefore there is a
vector bundle $V_0$ on $A_0$ such that
$$X_0=\bP(V_0).$$
Now
$$
-K_{X_{0}} = \cO_{\bP(V_0)}(2)\otimes\alpha^* (\det V_0^*),
$$
hence
$$W_0= W \vert A_0= S^2 V_0 \otimes \det V_0^*.$$
Therefore we have $c_1(W_0)=0$, hence $c_1(W)=0$, $E$ being finite, and
this proves (4.8.2).
\smallskip \noindent
By the Leray spectral sequence we get
$$\chi(-K_X) = \chi( \alpha_*(-K_X)) - \chi (R^1 \alpha_* (-K_X)) +
\chi(R^2 \alpha_* (-K_X)). \leqno{(4.8.3)}$$
We next claim that
$$R^2 \alpha_* (-K_X) = 0. \leqno{(4.8.4)} $$
Of course, $R^2 \alpha_* (-K_X) = 0$ generically. Now let
$a \in A$ and assume that
$\dim\alpha^{-1} (a) = 2$. Let $F= \alpha^{-1} (a)$, equipped with
its reduced structure. Then
$$R^2 \alpha_* (-K_X)_a = 0$$
if
$$H^2(F, -K_X \vert F \otimes N_F^{* \mu}) = 0 \leqno(4.8.5)$$
for all $\mu \ge 0$. Let $S \subset F$ be an irreducible 2-dimensional
component. Then (4.8.5) comes down to show
$$H^2(S, -K_X \vert S \otimes N_F^{* \mu} \vert S) = 0. \leqno{(4.8.6)}$$
By Serre duality, this means that
$$
H^0 (S, 2K_X \vert S \otimes N_F^{\mu} \vert S \otimes N_S)
=0. \leqno{(4.8.7)}
$$
If $X_a = \varphi^{-1}(a) $ denotes the full complex-analytic fiber
(with natural structure), then $N^*_{X_a}$ is generated by global
sections. It follows that $N_F^{* {\mu_1}} \vert S$ and $N_S^{*
{\mu_2}}$ have non-zero sections for suitable $\mu_1, \mu_2 > 0$. If
therefore (4.8.7) does not hold, we conclude -- having in mind that $-K_X
\vert S$ is nef -- that $K_X \vert S \equiv 0$, $N_F \vert S \equiv 0$
and $N_S \equiv 0$. Then however any section of $N_F^{* {\mu_1}} \vert
S$ resp $N_S^{*{ \mu_2}}$ is free of zeroes which implies $\alpha^{-1}
(a)= S$ set-theoretically. Then $N_S^*$ clearly cannot be numerically
trivial, since $N^*_{X_a}$ is generated by at least two sections and
$N^*_{X_a} \vert S = N_S^{* \lambda} \vert S$ generically suitable for
$\lambda$. Hence (4.8.7) holds and (4.8.4) is proved.

In completely the same way we prove that
$$R^2 {\alpha_*} (\cO_X) = 0. \leqno{(4.8.8)}$$
Since -- as already seen -- $\alpha$ is generically a $\bP_1$-bundle
with at most finitely many singular fibers, $R^1 \alpha_* (\cO_X)$ is
a torsion sheaf. Together with (4.8.8) and the Leray spectral sequence we
deduce
$$
\dim H^2(X,\cO_X) = 1. \leqno{(4.8.9)}
$$
Since $\dim H^1 (X, \cO_X) = q(X) = 2$, we obtain from (4.8.9):
$$
\chi(X, \cO_X)=0. \leqno{(4.8.10)}
$$
Since $-K_X$ is nef, we have $(-K_X)^3 \ge 0$. If however $(-K_X)^3
> 0$, then the holomorphic Morse inequalities imply the projectivity
of $X$ so that $(-K_X )^3 = 0$. Riemann-Roch and (4.8.10) therefore yield
$$
\chi(-K_X) = 0. \leqno{(4.8.11)}
$$
Then (4.8.3), (4.8.4) and (4.8.11) imply
$$
\chi(\alpha_*(-K_X)) \ge 0 \leqno{(4.8.12)}
$$
and $\chi( \alpha_*(-K_X)) =0$ if and only if $R^1 \alpha_* (-K_X) = 0.$

In order to bring $W$ into the game via (4.8.1), we first show:
$$
\chi(W) \le 0. \leqno{(4.8.13)}
$$
In fact, Riemann-Roch and (4.8.2) say
$$
\chi(W) = -c_2 (W)
$$
Suppose $\chi(W) >0$. Then $c_2 (W) < 0$.
Therefore the Bogomolov inequality
$$
c_1^2(W) \le 4 \, c_2 (W)
$$
is violated and $W$ is not semi-stable with respect to a fixed
K\"ahler metric $\omega$. Let $\cF \subset W$ be a maximal destabilizing
subsheaf w.r.t $\omega$; we may assume $\cF$ locally free of rank $1$ or $2$.
So
$$c_1 (\cF)\cdot\omega > 0.$$ First suppose that $\cF$ has rank $2.$
Then we have an exact sequence
$$ 0 \la \cF \la W \la {\cal I}_Z \otimes L \la 0$$
with a finite set $Z$ of length $l(Z)$ and a line bundle $L.$
Thus $ c_2(W) = c_2(\cF) + l(Z) - c_1(\cF)^2$ and $c_2(W) < 0$
yields
$$ c_1(\cF)^2 > c_2(\cF).$$
Now $c_1(\cF)^2 \leq 0$ since $A$ is non-algebraic and moreover
$c_1(\cF)^2 \leq 4c_2(\cF)$ since $\cF$ is $\omega-$semi-stable.
These last three inequalities give a contradiction.
The calculations in case that $\cF$ has rank $1$ are the same, working
with $W^*$ instead of $W.$ So (4.8.3) is verified.
\smallskip \noindent

Observe now that (4.8.12),(4.8.13) and
(4.8.1) imply
$$
\chi(\alpha_*(-K_X)) = \chi(W)=0$$ and $Q=0$, i.e.\
$\alpha_* (-K_X)=W.$ Thus $c_2(W)=0$, $W$ is semi-stable but not stable
and one can find $\cF \subset W$ as above with $c_1(\cF)=0$.
If $\rk \cF =1$, we have
$$
H^0(\alpha_*(-K_X) \otimes \cF^*) \not= 0
$$
proving (4.8.1).

In case $\rk \cF=2$ we consider the sequence $$0 \to \cF
\to W \to \cR \to 0. \leqno(4.8.14)$$
$\cR$ is a torsion
free coherent sheaf whose singular locus $Z=\Sing \cR$ is at
most finite (since $A$ has no compact positive dimensional
subvarieties). Consequently $\cR=\cJ \otimes L$ with an
ideal sheaf $\cJ$ such that Supp $(\cO_A/\cJ)=Z.$ With the
same argument as for $W$, we have
$$
c_2(\cF) \ge 0.
$$
The exact sequence $(4.8.14)$ yields
$$
c_2(W) = c_2(\cR) + c_2(\cF).
$$
Now $c_2(\cR)= \sharp Z$, counted with multiplicities. Therefore
$Z= \emptyset, \cR=L$ and $c_2(\cF)=0$.
Let $\zeta \in H^1(\cF \otimes L^*)$ be the extension class of
$$ 0 \la \cF \la W \la L \la 0. $$
If $\zeta \ne 0$, then
$H^1(\cF \otimes L^*) \ne 0,$ hence $H^0(\cF \otimes L^*) \ne 0$ or
$H^0(\cF \otimes \det \cF^* \otimes L) \ne 0$ by Riemann-Roch and duality.
Hence $H^0(W \otimes G) \ne 0$ for some $G$.
If $\zeta = 0,$ then $W = \cF \oplus L$, hence again $H^0(W \otimes L^*)
\ne 0$. \qed
\medskip

\noindent
Combining all our results we can state:

\claim 4.9.\ Theorem|Let $X$ be a compact K\"ahler 3-fold with $-K_X$
nef. Then the Albanese map is a surjective submersion. Moreover
$\pi_1(X)$ is almost abelian.
\endclaim

\proof. It remains only to treat the case $q = 0.$ By the results in
case $q > 0$, we may assume that
$$\tilde{q} (X)= \sup \{ q(\tilde{X})\,\vert~\tilde{X} \to X ~~
\hbox{finite \'etale}~\}=0.$$
Hence $H_1(X,\bZ)$ is finite, thus, passing to a finite \'etale cover,
we may assume\break $H_1(X,\bZ) = 0$. Denoting $G = \pi_1(X),$ we obtain
$$ G = (G,G).$$
By Paun [Pau96], see (4.20(b)), $G$ has polynomial growth. Hence
$G$ is almost nilpotent by Gromov's theorem, hence nilpotent.
Denoting $(C^k(G))$ the lower central series, we have
$C^2(G) = (G,G) = G, $ hence inductively $C^m(G) = G$ for all $m.$
But $C^n(G) = {e} $ for some $n$, $G$ being nilpotent. Hence $G =
{e}.$

\claim 4.10.\ Remark|{\rm It remains to investigate the structure of
compact
K\"ahler 3-folds $X$ with $-K_X$ nef and $q(X)=0$. Since our theory is
only up to finite \'etale covers we should assume
$$\tilde{q} (X)= \sup \{ q(\tilde{X})\,\vert~\tilde{X} \to X ~~
\hbox{finite \'etale}~\}=0.$$
In case $-K_X$ is hermitian semi-positive, the structure of $X$ is as
follows ([DPS96b]):} Up to a finite \'etale cover, $X$ is of one
of the following types:
\smallskip
\item{\rm a)} a Calabi-Yau manifold
\smallskip
\item{\rm b)} $\bP_1 \times S$, where $S$ is a K3 surface
\smallskip
\item{\rm c)} a rationally connected manifold.
\endclaim

In class c) we have all Fano $3$-folds, $3$-folds with $-K_X$ big and nef
but also e.g.\ $X=\bP_1 \times Y$, where $Y$ is $\bP_2$ blown up in 9
general points in such a way that $-K_X$ is hermitian semi-positive,
in particular $K_X^3=0$.
In case $-K_X$ nef one should expect the same result. Especially, if
$X$ is not projective then we should have $X=\bP_1 \times S$ with $S$
a non-projective $K 3$ surface.
Observe that rationally connected compact K\"ahler manifolds are
automatically projective ([Cam81]).\qed
\medskip

\noindent
We finish this section by proving a weak ``rational'' version of
our expectation in case $-K_X$ is nef and $X$ is projective.

\claim 4.11.\ Proposition|Let $X$ be a compact K\"ahler $3$-fold with
$-K_X$ nef and $K_X \not\equiv 0$. Assume $\tilde{q} (X)=0$. Then
either $X$ is rationally connected $($in particular
$\pi_1(X)=0)$ or there is a dominant meromorphic
map
$$f:X \rightharpoonup S$$
{\it to surface birational to a K3 or
Enriques surface $S$ with general fiber a rational curve (hence
$\pi_1(X)=0$ or $\pi_1(X)=\bZ_2) .$}
\endclaim

\proof. Assume $X$ not to be rationally connected. Since $-K_X$ is nef
and $K_X \not\equiv 0$, we have $\kappa (X)=- \infty$. Therefore
$X$ is uniruled. By [Cam92], [KoMM92] there is a
meromorphic fibration
$$f:X \rightharpoonup S$$
contracting the general rational curves of a given
(and fixed) covering family $(C_t)_{t \in T}$ of rational curves.
Properly speaking, we have $f(x) = f(y)$ for general points
$x,y \in X$, if $x$ and $y$ can be joined by a chain of rational
curves of type $C_t$. Since $X$ is not rationally connected, we have
dim $S>0.$
Next notice that dim $S=2$. In fact, if dim $S=1,$ then $S \simeq
\bP_1$ by $\tilde{q}(X)=0.$ Since the fibers of $f$ are rationally
connected, [KoMM92] gives rational connectedness.
Hence dim $S=2.$ Of course we may assume $S$ smooth.
We first verify that $\kappa(S) \le 0.$
\medskip

\noindent
Assume to the contrary that $\kappa(S) \ge 1.$

\noindent
We can obtain $f$ -- after possibly changing the family $(C_t)$ and
changing $S$ birationally and admitting rational double points on $S$
-- by a composite of birational Mori contractions and flips, say $X
\rightharpoonup X'$, and a Mori fibration
$$f':X'\to S$$
(just perform the Mori program on $X$).
Now $-K_{X'}$ has the following property (cf.\ [PS97]):
$$-K_{X'} \cdot C'\ge 0$$
for all curves $C'$ but a finite number of
rational curves. In particular $-K_{X'}$ is almost nef with non-nef
locus not projecting onto $S$, hence (4.7) applies and gives $\kappa(S)
\le 0.$
\medskip
In total we know $\kappa(S) \le 0.$ Since $\tilde{q}(S)=0,$
$S$ is birationally a K3, Enriques or rational surface. In the last case
$X$ is clearly rationally connected. \qed

\claim 4.12.\ Remark|{\rm (1) It seems rather plausible that methods
similar to those of [PS97] will prove that, in case $X$ is a
projective 3-fold with $-K_X$ nef, $\kappa(S)= 0$, $\tilde{q}
(X)=0$ and $X$ not rationally connected, the meromorphic map $f:X
\rightharpoonup S$ is actually a holomorphic map with $S$ \ul{being}\ a K3 or
Enriques surface and actually $f$ is a submersion, i.e.\ a $\bP_
1$-bundle. Then we see immediately that $X \simeq \bP_1 \times S$.  Of
course one difficulty arises from the fact that $S$ contains some
rational curves.  In the non-algebraic case one would further need a
more complete ``analytic Mori theory''.
\smallskip \noindent (2) We
discuss the K\"ahler analogue of (4.21). So let $X$ be a compact
K\"ahler threefold with $-K_X$ nef and $K_X \not \equiv 0.$ This first
difficulty is that $X$ might be simple, i.e.\ there is no positive
dimensional proper subvariety through the general point of $X.$ Ruling
out this potential case (which is expected not to exist), we conclude
by [CP00] that $X$ is uniruled. As in (4.21) we can form the rational
quotient $f: X \rightharpoonup S.$ Again $\dim S \ne 1$, because
otherwise $S = \bP_1$ and clearly $X$ is projective ($X$ cannot carry
a 2-form).  So $S$ is a non-projective surface and automatically
$\kappa (S) \leq 1.$ In fact, otherwise $X$ is algebraically connected
(any two points can be joined by a chain of compact curves) and
therefore projective by [Cam81]. One would like to exclude the case
$\kappa (S) = 1.$ However we do not know how to do this at the moment.
The method of 4.7 does not work because we do not have enough curves
in $S.$ What we can say is the following. $S$ admits an elliptic
fibration $h:S \to C \simeq \bP_1.$ Since $S$ is not algebraic, $h$
has no multi-section.  Now $f$ is almost holomorphic ([Cam92]), i.e.\
there are $U \subset X$ and $V \subset S$ Zariski open such that $f:U
\to V$ is holomorphic and proper. Let $F$ be a general fiber of $h$.
Let $A=S \ssm V.$ Then $A \cap F = \emptyset$ because $h$ has no
multi-sections, therefore $f$ is holomorphic over $F$ and the
composition yields a holomorphic map $X \la C.$ Let $X_F=f^{-1}(F).$
Then $X_F$ is a ruled surface over an elliptic curve of the form
$\bP(\cO \oplus L)$ with $L$ of degree $0$ but not torsion.}
\endclaim

In the second part of this section we investigate more generally the
structure of normal projective
varie\-ties $X$ such that $-K_X$ is pseudo-effective. Our leitfaden is
the following

\claim 4.13.\ Problem|Let $X$ and $Y$ be normal projective
$\bQ$-Gorenstein varieties. Let $\varphi: X \to Y$ be a
surjective morphism. If $-K_X$ is pseudo-effective $($almost nef$)$, is
$-K_Y$ pseudo-effective $($almost nef$)\,$?
\endclaim

\noindent
This problem in general has a negative answer:

\claim 4.14.\ Example|Let $C$ be any curve of genus $g \ge 2$. Let $L$ be
a line bundle on $C$ and put $X=\bP(\cO \oplus L)$.
Then we have
$$\eqalign{
H^0(-K_X)&=H^0 (X,\cO_{\bP(\cO \oplus L)}(2)
\otimes \pi^* (L^* \otimes -K_C))\cr
&=H^0(C, S^2(\cO \oplus L) \otimes
L^* \otimes -K_C) \not= 0,\cr}
$$
if deg $L^{*} > 3g - 2$.  So $-K_X$ is effective, hence
pseudo-effective (hence also almost nef by (1.5)), however $-K_C$ is
not pseudo-effective.
\endclaim

It is easily checked in this example that $-K_X \cdot B \ge 0$ for all
curves $B \subset X$ with the only exception $B=C_0$, $C_0$ the unique
exceptional section of $X$.  In particular the non-nef locus of $-K_X$
projects onto $C$. This leads us to reformulate (4.3) in the following way:

\claim 4.3.a.\ Problem|Assume in (4.3) additionally that the non-nef
locus of $-K_X$ does not project onto $Y$. Is $-K_Y$ pseudo-effective?
\endclaim

In case $-K_X$ is pseudo-effective the answer to (4.3.a) is positive at
least in the case
of submersions but with a slightly stronger assumption as in (4.3.a),
replacing the
non-nef locus by the zero locus of the multiplier ideal sheaf associated
with a metric
of minimal singularities (4.15 below) while in the almost nef case we have
onlya weak answer
(4.17 below), but for general $\phi$ and dealing with the non-nef locus.

\claim 4.15. Theorem|Let $X$ and $Y$ be compact K\"ahler
manifolds. Let $\varphi: X \to Y$ be a surjective submersion. Suppose that
$-K_X$
is pseudo-effective and that the zero locus of the multiplier ideal of a
minimal metric of
$-K_X$ does not
project onto $Y.$ Then $-K_Y$ is pseudo-effective.  \endclaim

\proof. This is a consequence of (2.6.1). In fact, apply (2.6.1) with
$L=-K_X$. \qed

\claim 4.16.\ Definition|A $(\bQ$-$)$line bundle $L$ on a normal
$n$-dimensional projective variety $X$ is {\it generically nef} if
$$L\cdot H_1 \cdots H_{n-1} \ge 0$$
for all ample divisors $H_i$ on $X$.
\endclaim

\claim 4.17.\ Theorem|Let $X$ and $Y$ be normal projective $\bQ-$Gorenstein
varieties and let  $\varphi:X \to Y$ be surjective. Let $-K_X$
be almost nef with non-nef locus $B$.  Assume $\varphi (B) \not= Y$.
Then $-K_Y$ is generically nef.
\endclaim

\proof. We will use the method of [Zh96] in
which Zhang proves the surjectivity of the Albanese map for projective
manifolds with $-K_X$ nef. In generalisation of Prop.~1 in [Zh96] we
claim

\noindent
(4.17.1) Let $\pi:X \to Z$ be a surjective morphism of smooth
projective varieties. Then there is no ample divisor $A$ on $Z$ such
that $-K_{X \vert Z} - \delta \varphi ^*(A)$ is pseudo-effective
(almost nef) with non-nef locus $B$ not projecting onto $Z$ for some
$\delta >0$, unless $\dim Z=0$.

\noindent
{\it Proof of $(4.17.1)$.} The proof is essentially the same as the one
of Theorem 2 in [Miy93], where Miyaoka proves that $-K_{X \vert Z}$
cannot be ample. Just replace $\Sing(\pi)$ in Lemma 10 by $\Sing(\pi)
\cup B$. Then, assuming the existence of $A$ and $\delta$, the old
arguments work also in our case.  This proves (4.17.1).

Now, coming to our previous situation, we argue as in [Zh96]. Take $C
\subset Y$ be a general complete intersection curve cut out by
$m_1H_1, \ldots, m_{n-1} H_ {n-1}$, with $H_i$ ample, $m_i >0$. Let
$X_C= \varphi^{-1} (C)$. By Bertini, $C$ and $X_C$ are smooth.
Applying (4.17.1) to $X_C \to C$ it follows that
$$\displaystyle -K_{X_C \vert C} - \delta \varphi^*(A)$$
is never pseudo-effective (almost nef) with non-nef locus not
projecting onto $C$ for any choice of $A$ and $\delta$. On the other
hand
$$\displaystyle -K_{X_C / C} = -K_{X/Z} \vert X_C = -K_X \vert X_C +
\varphi^* (K_Z \vert C).$$
Since $-K_X$ is pseudo-effective with non-nef locus not projecting
onto $Z$, so does $-K_X \vert X_C$ and we conclude that $K_Z \vert C$
cannot be ample, i.e.\ $K_Z \cdot C \le 0$, which was to be proved. \qed

\claim 4.18.\ Corollary|Let $X$ be a normal projective
$\bQ$-Gorenstein variety. Assume $-K_X$ almost
nef with non-nef locus $B$.
\smallskip
\item{\rm a)} If $\varphi:X \to Y$ is a surjective morphism to a
  normal projective $\bQ$-Gorenstein variety $Y$ with
  $\varphi(B) \not= Y$, then $\kappa (Y) \le 0.$

\smallskip
\item{\rm b)} The Albanese map $\alpha:X \to \Alb(X)$ is surjective,
  if $\alpha (B) \not= \alpha (X).$ \vskip0pt
\endclaim

The next question asks for the fundamental group of varieties $X$ with
$-K_X$ pseudo-effective and small non-nef locus. Here of course we
need more assumptions on the singularities. A {\it terminal
  (canonical)} $n$-fold is a normal projective $\bQ$-factorial
variety with at most {\it terminal (canonical)} singularities.

\claim 4.19.\ Question|{\rm Let $X$ be a terminal (or canonical) variety with
$-K_X$ pseudo-effective (almost nef). Assume that $\alpha (B) \not=
\alpha(X),$ with $\alpha: X \to \Alb (X)$ the Albanese and $B$ the
non-nef locus of $-K_X$. Is $\pi_1 (X)$ almost abelian (i.e.\ abelian
up to finite index)?}
\endclaim

\noindent
If $\alpha:X\to Y$ is a morphism, we will say that a line bundle $L$
on $X$ is {\it properly pseudo-effective (properly almost
nef)} with respect to $\alpha$, if $L$ is pseudo-effective (resp.\ nef)
and the non-nef locus satisfies $\alpha (B) \not= \alpha (X)$.

\claim 4.20.\ Remarks|{\rm
\smallskip
\item{\rm a)} If $-K_X$ is hermitian semi-positive, then $\pi_1(X)$ is
  almost abelian by [DPS96b].

\smallskip
\item{\rm b)} If $-K_X$ nef, then $\pi_1(X)$ has at most polynomial
  growth by Paun [Pau97].

\smallskip
\item{\rm c)} If $X$ is a projective surface such that $-K_X$ is properly
  pseudo-effective with respect to the Albanese map, then $\pi_1(X)$
  is almost abelian. This is a consequence of the Kodaira-Enriques
  classification and is proved as follows. First it is clear that
  $\kappa(X) \le 0$. If $\kappa(X)=0$, there is nothing to proved, $X$
  being birationally a torus, hyperelliptic, K3 or Enriques. So let
  $\kappa(X)=-{\infty}$.  By (4.4),(4.5), $X$ does not admit a map to
  a curve of genus $\ge2$, so that $X$ is either rational or its
  minimal model is a ruled surface over an elliptic curve. Thus
  $\pi_1(X)=0$ or $\pi_1(X)=\bZ^2$.\qed}  \endclaim

\noindent
In dimension 3, (4.20) has still a positive answer, at least if $X$ is
projective.

\claim 4.21.\ Theorem|Let $X$ be a terminal 3-fold such that $-K_X$ is
properly pseudo-effective with respect to the Albanese map. Then
$\pi_1(X)$ is almost abelian.
\endclaim

\proof. By [Mor88] there exists a finite sequence $\varphi:X
\rightharpoonup X'$ of extremal birational "divisorial" contractions
and flips such that either $K_{X'}$ is nef or $X'$ carries a Fano
fibration $\varphi:X' \to Y$, i.e.\ an extremal contraction with $\dim
Y < \dim X'$. Since extremal contractions and flips leave the
fundamental group unchanged, we have
$$\displaystyle \pi_1 (X)= \pi_1(X').$$

\noindent
a) We claim that $-K_{X'}$ is properly pseudo-effective with respect
to Albanese.  For that we need to prove the following. If $\lambda: X
\to Z$ is an extremal divisorial contraction (i.e.\ $\lambda$
contracts a divisor) or if $\lambda:X \rightharpoonup X^+$ is a flip,
then $-K_Z$ (resp $-K_{X^+}$) is properly pseudo-effective if $-K_X$
is properly pseudo-effective.  Indeed, if $\lambda:X \to Z$ is
divisorial, then $\lambda_{\ast}(-K_X)=-K_Z$ as $\bQ$-Cartier divisors
and moreover $\lambda_{\ast}:N^{1}(X) \to N^{1}(Z)$ maps effective
divisors to effective divisors, hence pseudo-effective divisors to
pseudo-effective divisors. If $\lambda:X \rightharpoonup X^+$ is
small, we still have a natural map
$$\lambda_{\ast}:N^1(X) \to N^1(Z)$$
with the same properties as
above. In fact, if $\cL \in \Pic(X)$, then let $C \subset X$
and $C^+ \subset X^+$ be the 1-dimensional indeterminacy sets so that
$X \ssm C \simeq X^+ \ssm C^+$, and consider $\cL^+=
\lambda_ {\ast}(\cL\vert X \ssm C)$. $\cL^+$ can
be extended to a reflexive sheaf on $X^+$, however since $X$ is
$\bQ$-factorial, so does $X^+$ ([KMM87]), and some
$(\cL^+)^m$ extends to a line bundle $\tilde \cL ^{[m]}$
on $X^+$. Now
$$H^0(X^+, \tilde \cL^{[m]})= H^0(X^+ \ssm
C^+, (\cL^+)^m)= H^0(X \ssm C, \cL^m)=
H^0(X, \cL^m).$$
In total $-K_{{X}^+}$ is again pseudo-effective.

It is now clear that if $-K_X$ is {\it properly} pseudo-effective
with respect to Albanese, then so is $-K_Z$ (resp.\ $-K_{X^+}$).
\medskip

\noindent
b) Now let $\varphi:X' \to Y$ be an extremal contraction with $\dim
Y \le 2$.

By a) we may assume $X=X'$.
If $\dim Y=0$, then $X$ is $\bQ$-Fano, hence rationally connected
by [KoMM92], and in particular $\pi_1(X)=0$.
If $\dim Y=1$, then $Y$ is a smooth curve with genus $g \le 1$ by
(4.17). Therefore $\pi_1(X)= \pi_1(Y)=0$ or $\bZ^2$.
If $\dim Y=2$, then $\kappa (\hat Y) \le 0,$ where $ \hat Y$ is a
desingularisation, again by (4.17). If $\kappa (\hat Y)=0$, then $\hat Y$ is
birational to a torus, a hyperelliptic surface, a $K 3$ surface or an
Enriques surface. Thus $\pi_1(\hat Y)$ is (almost) abelian. Since $Y$
has at most quotient singularities, we have $\pi_1(\hat Y)= \pi_1(Y)$,
hence $\displaystyle \pi_1(X)= \pi_1(Y)$ is (almost) abelian.

\noindent
c) Finally assume $K_X$ nef. Since $-K_X$ is pseudo-effective,
we must have $K_X \equiv 0$. Therefore $mK_X = \cO_X$ for suitable
$m>0$, and $\pi_1 (X)$ is almost abelian by [Kol95]. \qed

\claim 4.22.\ Remarks|{\rm
\smallskip
\item{\rm a)} If $-K_X$ is properly almost nef in (4.18) instead of
  properly pseudo-effective, then our conclusion still holds. The
  proof is essentially the same.  The only change concerns the
  invariance under flips. This follows from [KMM87, 5-1-11].
\smallskip
\item{\rm b)}
  To prove (4.21) in any dimension with the methods presented here,
  would require several deep things. First of all we would need the
  minimal model program working in any dimension.
\item{}
  Second we need to know that, given an extremal contraction $\varphi:
  X \to Y$ with $<\dim X $, then $-K_Y$ is properly pseudo-effective,
  if $-K_X$ is properly pseudo-effective; i.e.\ we would need a
  positive answer to (4.13a)) at least in the case of an extremal
  contraction.
\item{}
  And last we would need to know that $\pi_1(X)$ is almost abelian if
  $K_X \equiv 0$. This is well known if $X$ is smooth [Bea83] but hard
  if $X$ is singular, $\dim X \ge 4$. Compare [Pet93].\qed
\vskip0pt
}
\endclaim

\section{\S5. Threefolds with pseudo-effective canonical classes}

If $X$ is a smooth projective threefold or more generally a normal
projective threefold with at most terminal singularities such that
$K_X$ is pseudo-effective, then actually some multiple $mK_X$ is
effective, i.e.\ $\kappa (X) \geq 0.$ This is one of the main results
of Mori theory and is in fact a combination of Mori's theorem that a
projective threefold $X$ has a model $X'$ with either $X'$ uniruled or
with $K_{X'}$ nef and of Miyaoka's theorem that threefolds with $K_X$
nef have $\kappa (X) \geq 0.$ In particular reduction to char $p$ is
used and it is very much open whether the analogous result holds in
the K\"ahler category. In this section we give a very partial result
in this direction.

\claim 5.1.\ Lemma| Let $X$ be a normal compact K\"ahler space with at
most isolated singu\-la\-rities and $L$ a pseudo-effective line bundle on
$X.$ Let h be a singular metric on $L$ with curvature $\Theta_h(L)
\geq 0$ (in the sense of currents).  Let $\varphi$ be the weight
function of~$h.$ Assume that the Lelong numbers satisfy
$\nu(\varphi,x) = 0$ for all $x \in X$ but a countable set.
Then $L$ is nef.
\endclaim

\proof. See [Dem92], Corollary 6.4. \qed
\medskip

Now consider a pseudo-effective line bundle $L$ which is not nef.
Choose a metric $h_0$ with minimal singularities in the sense of
Theorem 1.5 and denote by $T_0$ its curvature current. Let
$$
E_c = \{x \in X \vert \nu(x,\varphi_0(x)) \geq c\}. $$
By [Siu74] (see also [Dem87]), $E_c$ is a closed analytic set.
Then Lemma 5.1 tells us that $\dim E_c \geq 1$ for sufficiently
small $c > 0.$ Furthermore we notice

\claim 5.2.\ Lemma|There exists $c>0$ such that $L \vert E_c $ is not nef.
Moreover, given an irreducible codimension $1$ (in $X$) component $D
\subset E_c,$
the line bundle $L \otimes \cO_X(-aD) \vert D$ is pseudo-effective
for a suitable $a > 0.$
\endclaim

\proof. In fact, by [Pau98b, Th\'eor\`eme 2], a closed positive
current $T$ has
a nef coholomogy class $\{T\}$ if and only if the restriction
of $\{T\}$ to all components $Z$ of all sets $E_c(T)$ is nef.
Thus, as $\{T_0\}$ is not nef, there must be some $E_c$ such that
$\{T_0\}|E_c$ is not nef. The second assertion follows from
Siu's decomposition
$$
T_0=\sum a_jD_j+R,
$$
where the $D_j$ are irreducible divisors, say $D_1=D$, and $R$ is a closed
positive current such that $\codim E_c(R)\ge 2$ for every $c$. Then
$R|D$ is pseudoeffective (as one sees by applying the main regularization
theorem for $(1,1)$-currents in [Dem92]), and $D_j|D$ is pseudo-effective
for $j\ge 2$, thus $\{T_0-a_1D_1\}|D$ is pseudoeffective.\qed

\claim 5.3.\ Corollary| Let $S$ be a compact K\"ahler surface and
$L$ a line bundle on $S.$ Assume that $L$ is pseudo-effective and that
$L\cdot C \geq 0$ for all curves $C \subset S.$ Then $L$ is nef.
\endclaim

\proof. Introduce a singular metric $h$ on $L$ whose curvature current
is positive (see 5.1). Let $c > 0$ and let $E_c$ be the associated Lelong
set of the weight function of $h.$ Then by 5.1 we can find $c$ such that
$\dim E_c = 1.$ Suppose that $L$ is not nef. Then by 5.2 there exists a
curve $C \subset E_c$ such that $L \cdot C < 0.$ This is a contradiction.
\qed

Of course this proof does not extend to dimension $3$, because now we
might find a surface $A$ such that $L\vert A$ is not nef. However
$L \vert A$ might not be pseudo-effective, so there is no conclusion.
On the other hand, the situation for $L = K_X$ is much better:

\claim 5.4.\ Theorem|Let $X$ be a $\bQ$-factorial normal
$3$-dimensional compact K\"ahler space with at most isolated
singularities. Suppose that $K_X$ is pseudo-effective but not nef.
Then there exists an irreducible curve $C \subset X$ with $K_X \cdot C
< 0.$ \smallskip \noindent If $X$ is smooth, then $C$ can be choosen
to be rational, and there exists a surjective holomorphic map $f: X
\longrightarrow Y$ with connected fibers contracting $C$ to a
point, such that $-K_X$ is
$f$-ample.  \endclaim

\proof. We only need to prove the existence of an irreducible curve
$C$ with $K_X \cdot C < 0;$ the rest in case $X$ is smooth follows
from [Pet98, Pet99].  \smallskip By Lemma 5.2 a suitable set $E_c$ contains
a positive-dimensional component $S$ with $K_X \vert S$ not nef.  If
$\dim S = 1,$ we are done, so suppose that $S$ is an irreducible
surface. Let $\mu: \tilde S \longrightarrow S$ denotes normalisation
followed by the minimal desingularisation. Let $L = \mu^*(K_X \vert
S)$. Then we need to show that there is a curve $C \subset \tilde S$
such that
$$L \cdot C < 0.$$
Since $L$ is not nef, this is clear if $\tilde S$ is projective. So
we may assume $\tilde S$ non-algebraic.
\medskip
By Lemma 5.2 we find a positive number $a$ such that $L + \mu^*(N_S^{*a})$
is pseudo-effective, hence by adjunction
$$ \mu^*(K_X^{(1+a)} \vert S - K_S^a) $$
is pseudo-effective. On the other hand by subadjunction,
$\mu^*(K_S) = K_{\tilde S} + B$
with $B$ effective. Since $K_{\tilde S}$ is effective, $\tilde S$ being
non-algebraic, we conclude that $L$ is pseudo-effective.
Now we apply 5.3 and conclude. \qed

The reason why we can construct a contraction in Theorem 5.3 is as
follows. If $X$ is a smooth K\"ahler threefold with $K_X \cdot C < 0$,
then $C$ moves in a positive-dimensional family, and therefore one can
pass to a non-splitting family.  This family provides the contraction,
see [CP97]. These arguments are likely to work also in the Gorenstein
case but break down in the presence of non-Gorenstein singularities.
Here it can happen that $C$ does not deform and new arguments are
needed.

\section{\S6. Pseudo-effective vector bundles}

In this section we discuss pseudo-effective vector bundles with
special emphasis on the tangent bundle of a projective manifold.

\claim 6.1.\ Definition|Let $X$ be a projective manifold and $E$ a holomorphic
vector bundle on $X.$ Then $E$ is {\it pseudo-effective}, if
$\cO_{\bP(E)}(1)$
is pseudo-effective and the union of all curves $C$ with
$\cO(1) \cdot C < 0$ $($i.e.\ the non-nef locus of the almost nef line
bundle $\cO(1))$ is contained in a countable union of subvarieties which do
not project onto $X$.
\endclaim

\claim 6.2.\ Remark| {\it  Notice that $E$ is pseudo-effective if and only if
$\cO_{\bP(E)}(1)$ is pseudo-effective and additionally there is a
countable union $S$ of proper subvarieties of $X$ such that $E \vert
C$ is nef for every curve not contained in $S.$}
\endclaim

\noindent
We have the following cohomological criterion for pseudo-effectivity.

\claim 6.3.\ Proposition|Let $X$ be a projective manifold and $E$ a
vector bundle on $X.$ Then $E$ is pseudo-effective if and only if
there exists an ample line bundle $A$ on $X$ and  positive integers
$m_0$ and $n_0(m)$   such that $H^0(X,S^n(S^mE \otimes A))$ generically
generates
$S^n(S^mE \otimes A)$ for all $m \ge m_0$ and $n \geq n_0(m)$.
\endclaim

\proof. If such an $A$ exists, then clearly $E$ is pseudo-effective,
using (1.2).  Working in the other direction, we choose an ample line
bundle $H$ on $\bP(E)$ by $H = \cO(1) \otimes
\pi^*(A)$ with $A$ ample on $X.$ Now let $x \in X$ be a very general
point and let $F = \pi^{-1}(x)$ the fiber over $\pi: \bP(E) \to
X.$ We must prove that
$$
H^0(\bP(E),(\cO(m) \otimes H)^{\otimes k} ) \to H^0(F,(\cO(m) \otimes
H)^{\otimes k}\vert F)$$
is surjective for $k \gg 0.$
Suppose first $E$ nef. Then reduce inductively to $\dim X = 1;$ the
necessary $H^1-$vanishing is provided
by Kodaira's vanishing theorem for sufficiently large $k.$ In case $\dim X
= 1,$ the ideal sheaf
${\cI}_F$ is locally free of rank $1$ and
$$H^1(\bP(E),{\cI}_F \otimes (\cO(m) \otimes H)^k) = 0 $$
again holds by Kodaira for large $k.$
In the general case one
introduces multiplier ideal sheaves associated with singular metrics on
$\cO(m)$ whose support do not meet $F$ and
substitutes Kodaira's vanishing theorem by Nadel's vanishing theorem.
We leave the easy details to the reader. \qed

\claim 6.4.\ Definition|Let $X$ be a projective manifold and $E$ a
vector bundle on $X.$ Then $E$ is said to be {\it almost nef}, iff
there is a countable family $A_i $ of proper subvarieties of $X$ such
that $E \vert C$ is nef for all $C \not \subset \bigcup _{i} A_i.$ The
non-nef locus of $E$ is the smallest countable union $S$ of analytic
subsets such that $E \vert C$ is nef for all $C \not \subset S.$
\endclaim

\noindent
It is immediately seen that $E$ is almost nef if and only if
$\cO_{\bP(E)}(1)$ is almost nef with non-nef locus not projecting onto
$X.$ Hence we have in analogy to (3.3) the following

\claim 6.5.\ Proposition| Let $X$ be a projective manifold and $E$ a
vector bundle on $X.$ If $E$ is pseudo-effective, then $E$ is almost
nef.
\endclaim

\claim 6.6.\ Problem|Is every almost nef vector bundle pseudo-effective$\,$?
\endclaim

\noindent
 Here are some basic properties of pseudo-effective and almost nef
vector bundles.

\claim 6.7.\ Theorem|Let $X$ be a projective manifold and $E$ a vector
bundle on $X.$
\smallskip
\item{\rm a)} If $E$ is pseudo-effective (almost nef) and $\Gamma^{a}$
  is any tensor representation, then $\Gamma^{a}E$ is again
  pseudo-effective (almost nef). In particular $S^mE$ and $\Lambda^qE$
  are pseudo-effective (almost nef).
\smallskip
\item{\rm b)} If $E$ is almost nef and if $s \in H^0(E^*)$ is a non-zero
section, then $s$ has no
zeroes at all.
\smallskip
\item{\rm c)}
  If either $E$ is pseudo-effective or if $E$ is almost nef with non-nef locus
  $S$ having codimension at least $2,$ and if $\det E^*$ is almost
  nef, then $E$ is numerically flat i.e.\ both $E$ and $E^*$ are nef,
  and then $E$ has a filtration by hermitian flat bundles.
\vskip0pt
\endclaim

\proof. a) By standard representation theory it is sufficient to prove
the statement for $E^{\otimes m}.$ But in that case the claim is immediate.

\noindent
b) Let $S \subset X$ be the non-nef locus. Suppose $s(x) = 0.$ Choose
a curve $C \not \subset S$ such that $x \in C.$ Then $E\vert C$ is
nef, on the other hand $E^* \vert C$ has a section with zeroes. So $s
\vert C = 0.$ Varying $C$ we conclude $s = 0.$

\noindent
c) First notice that our claim is easy for line bundles $L,$ even
without assumption on the codimension: if both $L$ and $L^*$ are
almost nef, then clearly $L \cdot H_1 \cdots H_{n-1}= 0$
for all ample line bundles $H_i$ on $X$ ($n = \dim X).$ Thus $L \equiv
0.$ In particular we conclude from (a) that $\det E \equiv 0.$
\medskip \noindent
We now treat the case when $E$ is almost nef with $\codim S \geq 2.$

\noindent
c.1) First we claim that $E$ is $H$-semistable for all ample line
bundles $H.$ In fact, otherwise we find some $H$ and a torsion free
subsheaf $\cS \subset E$ such that $c_1(\cS) \cdot H^{n-1} > 0.$ Now,
assuming $H$ very ample, let $C$ be a general complete intersection
curve cut out by $H.$ Then $\cS \vert C$ is locally free and $c_1(\cS
\vert C) > 0.$ On the other hand $E \vert C$ is numerically flat
(since it is nef and $\det E^*\vert C$ is nef,
see [DPS94]). This is impossible. So $E$ is $H$-semistable for all $H.$

\noindent
c.2) As a consequence we obtain the inequality
$$ 0 = (r-1)c_1^2(E) \cdot H^{n-2} \leq 2rc_2(E) \cdot H^{n-2}$$
for all $H$ ample, $r$ denoting the rank of $E.$

\noindent
c.3) Next suppose that $E$ is $H$-stable for some $H$ (and still that
$\codim S \geq 2).$ Let $Y$ be a general surface cut out by hyperplane
sections in $H$ (again assume $H$ very ample). Then $Y \cap S =
\emptyset$ by our assumption on the codimension of $S,$ hence $E \vert
Y$ is nef, hence numerically flat and we conclude
$$c_2(E) \cdot H^{n-2} = 0. \leqno (6.7.1) $$
Since $E$ is $H$-stable, $E$
is Hermite-Einstein and from (6.7.1) we deduce (see e.g.\ [Kob87, p.115])
that $E$ is numerically flat.  Hence we may assume that $E$ is
$H$-semistable for all $H$ but never $H$-stable.

\noindent
c.4) Fix some ample line bundle $H$ and let $\cS$ the maximal
$H$-destabilizing subsheaf of $E,$ so that
$$c_1(\cS) \cdot H^{n-1}= 0. \leqno (6.7.2)$$
On the other hand the generically surjective map $E^* \to \cS^*$ proves that
$$
c_1(\cS^*) \cdot H_1 \cdots H_{n-1} \geq 0
$$
for all ample divisors $H_i.$ Together with (6.7.2) this yields
$$c_1(\cS) = 0.$$
Now we follows the arguments of [DPS94].  Let $p =
{\rm rk}\cS.$ Then $\det \cS$ is a numerically flat line bundle, moreover
it is a subsheaf of $\Lambda^pE$, hence by (6.7~a), $\det \cS$ is a
subbundle of $\Lambda^pE,$ and thus by [DPS94,1.20], $\cS$ is a
subbundle of $E.$ Now $\cS$ being almost nef with trivial determinant,
an induction on the rank of $E$ yields the numerical flatness of $\cS.$
For the same reason the quotient bundle
$E/\cS$ is numerically flat, too, so that $E$ is numerically flat.

\noindent
c.5) Finally assume that $E$ is pseudo-effective and suppose $S \ne
\emptyset$. Choose a general smooth curve $C$ meeting $S$ in a finite
set.  Then $E_C$ is nef with $\det E_C \equiv 0.$ Therefore $E_C$ is
numerically flat. Let $x \in C \cap S$ with maximal ideal sheaf $\cI_x
\subset \cO_C$ and let $A$ be ample on $X$. Then, applying 6.3 and
having in mind that $S^mE \vert C'$ is not nef (with $C'$ a suitable
curve in $S$ passing through $x$), the image of the restriction map
$$
H^0(X,S^k(S^mE \otimes A)) \longrightarrow
H^0(C,S^k(S^m(E_C \otimes A))
$$
has non-zero intersection with
$$
H^0(C,S^k(S^m(E_C \otimes \cI_x) \otimes A))
$$
for $m$ large and sufficiently divisible and for large $k.$
This contradicts the numerical flatness of $E_C.$  \qed

\claim 6.8.\ Remark|{\rm Of course we expect that every almost nef
bundle $E$ with $\det E \equiv 0$ is numerically flat.  The above
considerations show that it is sufficient to prove this only on
surfaces and for bundles which are stable (for all polarisations).
We formulate the problem precisely as follows.}
\endclaim

\claim 6.9.\ Problem|Let $Y$ be a smooth projective surface and $E$ a
vector bundle on $Y$ of rank at least~$2$. Let $E$ be an almost nef vector
bundle on $X.$
If
$\det E \equiv 0$, show that $E$ is nef, hence numerically flat.
\endclaim

\noindent
It would be sufficient to prove that $c_2(E) = 0$ and also the proof of (6.7)
shows that one may assume $E$ to be $H-$stable for any ample polarisation $H$ on
$Y.$ It is a priori clear that $E$ is always $H-$semi-stable.

\noindent We now want to study projective manifolds with almost nef tangent
bundles. A class of examples is provided by the almost homogeneous
manifolds $X$, i.e.\ the automorphim group acts with an open orbit, or
equivalently, the tangent bundle is generically generated. A question
we have in mind is how far the converse is from being true.

\claim 6.10.\ Proposition|Let $X$ be a projective manifold with $T_X$
almost nef. Then \hbox{$\kappa (X)\!\le\!0.$}
If $\kappa (X)=0,$ then $K_X\equiv 0.$
\endclaim

\proof. Since $T_X$ is almost nef, $-K_X$ is almost nef, too. If
$\kappa (X) \ge 0,$ then $mK_X$ is effective, therefore $K_X \cdot H_1
\cdots H_{n-1} =0$ for all hyperplane sections $H_i$, hence $K_X
\equiv 0$. \qed

\claim 6.11.\ Proposition|Let $X$ be a projective manifold with $K_X
\equiv 0$ and $T_X$ pseudo-effective. Then $X$ is abelian after a suitable
finite
\'etale cover.
\endclaim

\proof.  By (6.7) $T_X$ is numerically flat  and the claim
follows by Yau's theorem. \qed

\claim 6.12.\ Corollary|Let $X$ be a Calabi-Yau or projective
symplectic manifold. Then neither $T_X$ nor $\Omega^1_X$ is
pseudo-effective.  Moreover the union of curves $C \subset X$ such
that $T_X \vert C$ (resp. $\Omega^1_X \vert C)$) is not nef is not
contained in a countable union of analytic sets of codimension at
least~$2$.
\endclaim

\noindent
Of course we expect (see Problem 6.6) that the tangent bundle of a Calabi-Yau or
symplectic manifold is not almost nef, but this seems rather delicate
already for K3 surfaces. Since both bundles are generically nef by
Miyaoka's theorem, we get examples of generically nef vector bundles (of
rank at least $2$) which are not pseudo-effective.

\claim 6.13.\ Example|{\rm Let $X \subset \bP_3$ be a general quartic (with
$\rho (X)=1).$ Then $T_X$ is not pseudo-effective by (6.11). More
precisely :
$$H^0 (X, S^k(S^{m} T_X \otimes \cO_X (1)))=0\qquad
\hbox{for all $m \ge 1$ and all $k \ge 1.$}
$$
Here $\cO_X (1)$ is the ample generator of $Pic (X) \cong \bZ.$
This is proved by direct calculation in [DPS96a] and actually shows that
$\cO_{\bP(T_X)}(1)$ is not pseudo-effective. The same argument works for
smooth hypersurfaces
$X \subset \bP_{n+1}$ of degree $n+2$ with $\rho (X)=1.$
\medskip Continuing this example we are now able to exhibit a generically
nef line bundle $L$ which
is not pseudo-effective. Namely, let
$$ L = \cO_{\bP(T_X)}(1).$$
As already noticed, $L$ is not pseudo-effective. So it remains to verify
that $L$ is generically nef which is
of course based on the generic nefness of $T_X.$ Let $H_1$ and $H_2$ be
ample divisors on $\bP(T_X),$ then we
must prove
$$ L\cdot H_1 \cdot H_2 \geq 0.$$
Let $m_i \gg 0$ and choose $D_i \in \vert m_i H_i \vert.$ Let
$$\tilde C = D_1 \cap D_2$$
and let $C = \pi(\tilde C)$, where $\pi: \bP(T_X) \to X$ is the projection.
Then certainly $\deg(\tilde C/C) \geq 2;$ moreover
$$ C \in \vert \cO_X(k) \vert,$$
where $k$ is so large that $T_X \vert C_0$ is nef for the general $C_0 \in
\vert \cO_X(1) \vert.$
Thus $\cO_{\bP(T_X \vert C_0)}(1) $ is nef.
By semi-continuity (for $\cO_{\bP(T_X \vert C_0)}(m) \otimes A$ with $A$
ample on $\bP(T_X)$), we
conclude that $\cO_{\bP(T_X \vert C)}(1)$ is at least pseudo-effective.
Since by reasons of degree, $\tilde C$ cannot be the exceptional section of
$\bP(T_X \vert C),$
we conclude that
$$ m_1 m_2 (L \cdot H_1 \cdot H_2) = c_1(\cO_{\bP(T_X \vert C)}(1)\vert
\tilde C ) \geq 0,$$
proving our claim.}\qed
\endclaim

\noindent We now study projective manifolds $X$ with $T_X$ almost nef and
$\kappa (X) = - \infty.$

\claim 6.14.\ Proposition|Let $X$ be a projective manifold with $T_X$
almost nef and $\kappa (X) = - \infty.$ Then $X$ is uniruled.
\endclaim

\proof. In fact, since $T_X$ is almost nef we have
$$
K_X \cdot H_1 \cdots H_{n-1}\leq 0$$
for all ample divisors $H_i$ on $X.$ We must have strict inequality
for some choice of $H_i$ because otherwise $K_X \equiv 0$ and $\kappa
(X) = 0.$ Now [MM86] gives the conclusion. \qed

\claim 6.15.\ Proposition|Let $X$ be a projective manifold with $T_X$
almost nef. Then
\smallskip
\item{\rm a)} the Albanese map is a surjective submersion.
\smallskip
\item{\rm b)} there is no surjective map onto a variety $Y$ with
$\kappa (Y) > 0.$
\endclaim

\proof. In fact, all holomorphic $1$-forms on $X$ have no zeroes,
and this proves a) (compare with [DPS94]). Point b) follows from the
analogous fact that sections in $S^m\Omega^p_X$ cannot have zeroes.\qed

\noindent For the following considerations we recall from (4.20) the
following notation:
$$
\tilde q(X) = \sup\{q(\tilde X)\,\vert~\tilde X \to X~\hbox{\rm
finite \'etale} \}.
$$

\claim 6.16.\ Proposition|Let $T_X$ be pseudo-effective and suppose that
$\pi_1(X)$ does not contain a non-abelian free subgroup. Then if
$H^0(\Omega^p_X) \ne 0$ for some $p,$ we have $q(\tilde X) > 0.$
\endclaim

\proof. The proof is contained in [DPS94, 3.10].\qed

\noindent Notice that $\pi_1(X)$ does not contain a free non-abelian
free subgroup if $-K_X$ is nef [DPS93].

\noindent
In order to make further progress we need informations on manifolds
with $T_X$ almost nef which do not have $p$-forms for all $p,$ the
same also being true for every finite \'etale cover. \qed

\claim 6.17.\ Conjecture|Let $T_X$ be almost nef. Assume that
$H^q(\tilde X,\cO_{\tilde X}) = 0$ for all finite \'etale covers and
all $q \geq 1.$ Then $\pi_1(X) = 0,$ in fact $X$ should even be
rationally connected.
\endclaim

\noindent
The evidence for the conjecture is that it holds if $T_X$ is nef by
[DPS94]; it is furthermore true if $X$ is almost homogeneous and it is
true in low dimensions, as we shall see below.  With the same
arguments as in [DPS94] we obtain

\claim 6.18.\ Proposition|Let $T_X$ be almost nef and suppose that
$\pi_1(X)$ does not contain a nonabelian free subgroup. If Conjecture
6.17 holds, then $\pi_1(X)$ is almost abelian, i.e.\ $\pi_1(X)$
contains a subgroup $\bZ^r$ of finite index.
\endclaim

\claim 6.19.\ Proposition|Let $X$ be a smooth projective surface with
$T_X$ almost nef and $\kappa (X) = - \infty.$
\smallskip
\item{\rm a)} If $q(X) > 0,$ then $X$ is ruled surface over an
  elliptic curve and every ruled surface over an elliptic curve has
  $T_X$ almost nef.
\smallskip
\item{\rm b)} The minimal model $Y$ has again $T_Y$ almost nef.
\smallskip
\item{\rm c)} Every rational ruled surface has an almost nef tangent bundle.
\endclaim

\proof. a) The first statement is clear from (6.15), the second part
is done as follows. Let $p: X \to C$ be the projection to the elliptic
curve. Then consider the tangent bundle sequence
$$ 0 \to T_{X/C} \to T_X \to \cO_X \to 0.$$
Then $T_{X/C} = -K_X,$ and the claim follows from the pseudo-effectivity
of $-K_X$, which is immediately checked by [Ha77,V.2].

\noindent
b) is obvious.

\noindent
c) Such an $X$ is actually almost homogeneous, see [Po69]. \qed

\noindent
We do not investigate the rather tedious problem of determining which
rational blow-ups have almost nef tangent bundles, and instead turn
ourselves to the case of $3$-folds.

\claim 6.20.\ Theorem|Let $X$ be a smooth projective $3$-fold with $T_X$
almost nef.  Then $\pi_1(X)$ is almost abelian and if $\tilde q(X) =
0,$ then $X$ is rationally connected unless (possibly)
$T_X$ is not pseudo-effective and $X$ is a $\bP_1-$bundle over a K3- or an
Enriques surface (this case should not exist). Suppose $\kappa (X) = -\infty.$
Then the finer structure of $X$ is as follows:
\smallskip
\item{\rm a)} if $q(X) = 2, $ then $X$ is a $\bP_1$-bundle over an
  abelian surface.
\smallskip
\item{\rm b)} if $q(X) = 1, $ then the Albanese map $\alpha : X \to A$
  is a fiber bundle over the elliptic curve with general fiber $F$
  having $T_F$ almost nef and $\kappa (F) = - \infty.$ \smallskip
  \itemitem{$\alpha)$} If $F$ is rational, then there is a factorisation
  $X \to Y \to A$ such that either $f: X \to Y$ is birational, in that
  case $f$ is a succession of blow-ups of \'etale multisections over
  $A.$ Or $\dim Y < \dim X,$ in that case $f: X \to Y$ and $g: Y \to
  A$ are both $\bP_1$-bundles or $\alpha = \phi$ is a $\bP_2$- or a
  $\bP_1 \times \bP_1$-bundle.
\smallskip
\itemitem{$\beta)$} If $F$ is a ruled surface over an elliptic curve, then
  $f$ is a $\bP_1$-bundle and $g: Y \to A$ is a hyperelliptic surface.
  Then there is an \'etale $2:1$-cover $\tilde X \to X$ such that
  $q(\tilde X) = 2.$
\vskip0pt
\endclaim

\proof. If $\kappa (X) = 0,$ then we have $K_X \equiv 0$ and $\pi_1(X)$
is abelian by [Bea83].
So we will now suppose $\kappa (X) = - \infty.$

\noindent
(I) First we treat the case $\tilde q(X) = 0.$ By (6.14) $X$ is
uniruled, so we can form the rational quotient $f: X \rightharpoonup Y
$ with respect to some covering family of rational curves [Cam81,
Cam92], [KoMM92]. Then $Y$ is a projective manifold with $\dim Y \leq
2.$ If $\dim Y = 0,$ then $X$ is rationally connected. If $\dim Y =
1,$ then $Y = \bP_1$ by $\tilde q(X) = 0,$ hence $X$ is rationally
connected by [KoMM92]. If finally $\dim Y = 2,$ then from $\tilde q(Y)
= 0$ and $\kappa (Y) \leq 0$ (6.15), we deduce that either $Y$ is rational
or a K3 resp.
an Enriques surface. In the rational case
it is easy to see and well-known
that $X$ is rationally connected.
So suppose $Y$ is not rational. Since an Enriques surface has a finite
\'etale cover which
is K3, we may assume that $Y$ is actually K3. Then $X$ carries a
holomorphic 2-form $\omega.$
Since $T_X$ is almost nef, $\omega$ cannot have zeroes (6.7)(b)).
Therefore $X$ is a
$\bP_1-$bundle over $Y$ by [CP00]. Now suppose $T_X$ pseudo-effective. Then
$f^*(T_Y)$ is
pseudo-effective as quotient of $T_X$. Since $c_1(f^*(T_Y)) = 0$, the
bundle $f^*(T_Y) $ is
numerically flat by (6.7), in particular $f^*(c_2(Y)) = 0$ which is absurd.

\noindent
(II) We will now assume $q(X) > 0$ and shall examine the structure of
the Albanese map $\alpha: X \to A.$ Once we have proved the structure
of $\alpha$ as described in 1) and 2), it is clear
that $\pi_1(X)$ is almost abelian. We already know that $\alpha $ is
a surjective submersion.

\noindent
a) If $\dim A = 2,$ then $\alpha$ is a $\bP_1$-bundle.

\noindent
b) Suppose now $\dim A = 2.$ Let $F$ be a general fiber of $\alpha.$
Then $\kappa (F) = - \infty,$ since $X$ is uniruled (6.14), moreover
$T_F$ is almost nef, hence either $F$ is rational or a ruled surface
over an elliptic curve (6.19). We will examine both cases by studying
a Mori contraction $\phi: X \to Z$ which induces a factorisation
$\beta: Z \to A$ such that
$$\alpha = \beta \circ \phi.$$
Note that all possible $\phi$ are classified in [Mor82].

\noindent
b.1) Assume that $F$ is rational. If $ \dim Z = 2, $ then $\phi$ is a
$\bP_1$-bundle and $\beta$ is again a $\bP_1$-bundle. If $\dim Z = 1,$
then $A = Z$ and $\alpha = \phi$ is a $\bP_2$- or a $\bP_1 \times
\bP_1$-bundle. If finally $\dim Z = 3,$ then $\phi$ must be the
blow-up of an \'etale multisection of $\beta$ and $T_Z$ is again
almost nef so that we can argue by induction
on $b_2(X).$

\noindent
b.2) If $F$ is irrational, then $\phi$ is necessarily a
$\bP_1$-bundle over $Z$ and then $\beta : Z \to A$ is an elliptic
bundle with $\kappa (Z) = 0,$ hence $Z$ is hyperelliptic. \qed

\section{References}

\bigskip

{\eightpoint

\bibitem [At57]&Atiyah, M.F.:& Vector bundles over elliptic curves;&
  Proc.\ London Math.\ Soc.\ {\bf 27} (1957) 414--452&

\bibitem [BPV]&Barth, W., Peters, C., Van de Ven, A.:& Compact Complex
  Surfaces;& Erg.\ der Math.\ 3.\ Folge, Band 4, Springer, Berlin
  (1984)&

\bibitem [Bea83]&Beauville, A.:& Vari\'et\'es k\"ahleriennes dont la
  premi\`ere classe de Chern est nulle;& J.\ Diff.\ Geom.\ {\bf 18}
  (1983) 745--782&

\bibitem [Cam81]&Campana, F.:& Cor\'eduction alg\'ebrique d'un espace
  analytique faiblement k\"ahl\'erien compact;& Inv.\ Math.\ {\bf 63}
  (1981) 187--223&

\bibitem [Cam92]&Campana, F.:& Connexit\'e rationnelle des
  vari\'et\'es de Fano;& Ann.\ Sci.\ ENS {\bf 25} (1992) 539--545&

\bibitem [CP97]&Campana, F.; Peternell, Th.:& Towards a Mori theory on
  compact K\"ahler $3$-folds, I;& Math.\ Nachr.\ {\bf 187} (1997) 29--59&

\bibitem [CP00]&Campana, F.; Peternell, Th.:& Holomorphic $2$-forms
  on complex $3$-folds;& J.\ Alg.\ Geom.\ {\bf 9} (2000) 223--264&

\bibitem [CPZ98]&Campana, F.; Peternell, Th.; Zhang, Q.:& On the
  Albanese maps of compact K\"ahler manifolds with nef anticanonical
  classes;& Preprint 1998&

\bibitem [De82]&Demailly, J.-P.:& Estimations $L^2$ pour l'op\'erateur
  $\overline\partial$ d'un fibr\'e vectoriel holomorphe semi-positif
  au-dessus d'une vari\'et\'e k\"ahl\'erienne compl\`ete;& Ann.\ Sci.\
  \'Ecole Norm.\ Sup.\ 4e S\'er. {\bf 15} (1982) 457--511&

\bibitem [De87]&Demailly, J.-P.:& Nombres de Lelong g\'en\'eralis\'es,
  th\'eor\`emes d'int\'egralit\'e et d'analyticit\'e;& Acta Math.\
  {\bf 159} (1987) 153--169&

\bibitem [De90]&Demailly, J.-P.:& Singular hermitian metrics on positive
  line bundles;& Proceedings Alg.\ Geom.\ Bayreuth 1990; Lecture Notes
  in Math.\ {\bf 1507}, Springer-Verlag (1992) 87--104&

\bibitem [De92]&Demailly, J.-P.:& Regularization of closed
  positive currents and intersection theory;& J.\ Alg.\ Geom.\ {\bf 1} (1992),
  361--409&

\bibitem [DPS93]&Demailly, J.-P.; Peternell, Th.; Schneider, M.:&
  K\"ahler manifolds with numerically effective Ricci-class;&  Comp.\
  Math.\ {\bf 89} (1993) 217--240&

\bibitem [DPS94]&Demailly, J.-P.; Peternell, Th.; Schneider, M.:&
  Compact complex manifolds with numerically effective tangent
  bundles;& J.\ Alg.\ Geom.\ {\bf 3} (1994) 295--345&

\bibitem [DPS96a]&Demailly, J.-P.; Peternell, Th.; Schneider, M.:&
  Holomorphic line bundles with partially vanishing cohomology;&
  Israel Math.\ Conf.\ Proc.\ vol.~{\bf 9} (1996) 165--198&

\bibitem [DPS96b]&Demailly, J.-P.; Peternell, Th.; Schneider, M.:&
  Compact K\"ahler manifolds with hermitian semi-positive anticanonical
  bundle;& Compositio Math.\ {\bf 101} (1996) 217--224&

\bibitem [EV92]&Esnault, H.; Viehweg, E.:& Lectures on Vanishing Theorems;&
  DMV Seminar, Band {\bf 20}, Birkh\"auser Verlag (1992)&

\bibitem [FiFo79]&Fischer, G.; Forster, O.:& Ein Endlichkeitssatz f\"ur
  Hyperfl\"achen aud kompakten komplexen R\"aumen;&
  J.\ f\"ur die reine und angew.\ Math.\ {\bf 306} (1979) 88--93&

\bibitem [GR84]&Grauert, H., Remmert, R.:& Coherent analytic
  sheaves;& Grundlehren der math.\ Wissenschaften {\bf 265}, Springer-Verlag,
  Berlin (1984)&

\bibitem [Ha77]&Hartshorne, R.:& Algebraic geometry;&
  Graduate Texts in Math.\ vol.\ 52 , Springer-Verlag 1977&

\bibitem [H\"or66]&H\"ormander, L.:& An introduction to Complex
  Analysis in several variables;& 1966, 3rd edition, North-Holland Math.\
  Libr., vol.7, Amsterdam, London (1990)&

\bibitem [KMM87]&Kawamata, Y.; Matsuki, K.; Matsuda, K.:&
  Introduction to the minimal model program;& Adv.\ Stud.\ Pure Math.\
  {\bf 10} (1987) 283--360&

\bibitem [Kob87]&Kobayashi, S.:& Differential geometry of complex vector
  bundles;& Princeton Univ.\ Press (1987)&

\bibitem [Kol86]&Koll\'ar, J.:& Higher direct images of dualizing sheaves I;&
  Ann.\ Math.\ {\bf 123} (1986) 11--42&

\bibitem [Kol95]&Koll\'ar, J.:& Shafarevich maps and plurigenera of
  algebraic varieties;& Princeton Univ.\ Press, 1995&

\bibitem [Kol96]&Koll\'ar,J.:& Rational curves on algebraic varieties;&
Springer, 1996&

\bibitem [KoMM92]&Koll\'ar, J.; Miyaoka, Y.; Mori, S.:& Rationally
  connected varieties;& J.\ Alg.\ Geom.\ {\bf 1} (1992) 429--448&

\bibitem [Kra75]&Krasnov, V.A.;& Compact complex manifolds without
  meromorphic functions (Russian);& Mat.\ Zametki {\bf 17} (1975)
  119--122&

\bibitem [LB92]&Lange, H.; Birkenhake,C.:& Complex abelian varieties;&
Springer, 1992&

\bibitem [Lel69]&Lelong, P.:& Plurisubharmonic functions and positive
  differential forms;& Gordon and Breach, New-York, and Dunod, Paris (1969)&

\bibitem [Miy93]&Miyaoka, Y.:& Relative deformations of morphisms and
  applications to fiber spaces;\break& Comm.\ Math.\ Univ.\ St.\ Pauli
  {\bf 42} (1993) 1--7&

\bibitem [Mor82]&Mori, S.:& Threefolds whose canonical bundles are not
  numerically effetive;& Ann.\ Math.\ {\bf 116} (1982) 133--176&

\bibitem [Mor88]&Mori, S.:& Flip theorem and existence of minimal
  models of $3$-folds;& J.\ Amer.\ Math.\ Soc.\ {\bf 1} (1988) 177--353&

\bibitem [MM86]&Miyaoka, Y; Mori, S.:& A numerical criteria for
  uniruledness;&  Ann.\ Math.\ {\bf 124} (1986) 65--69&

\bibitem [Mou97]&Mourougane, Ch.:& Images directes de fibr\'es en droites
  adjoints;& Publ.\ Res.\ Inst.\ Math.\ Sci.\ {\bf 33} (1997) 893--916&

\bibitem [Mou99]&Mourougane, Ch.:& Th\'eor\`emes d'annulation
g\'en\'eriques pour
  les fibr\'es vectoriels semi-n\'egatifs;& Bull.\ Soc.\ Math.\ Fr.\
  {\bf 127} (1999) 115--133&

\bibitem[Nad89]&Nadel, A.M.:& Multiplier ideal sheaves and
  K\"ahler-Einstein metrics of positive scalar curvature;& Proc.\ Nat.\
  Acad.\ Sci.\ U.S.A.\ {\bf 86} (1989) 7299--7300~~~and\hfil\break
  Annals of Math., {\bf 132} (1990) 549--596&

\bibitem[Nar66]&Narasimhan, R.:& Introduction to the theory of analytic
  spaces;& Lecture Notes in Math.\ no~{\bf 25}, Springer-Verlag (1966)&

\bibitem[NS95]&Namikawa, Y.; Steenbrink, J.H.M.:& Global smoothing of
  Calabi-Yau $3$-folds;& Invent.\ Math.\ {\bf 122} (1995) 403--419&

\bibitem[OT87]&Ohsawa, T.; Takegoshi, K.:& On the extension
  of $L^2$ holomorphic functions;& Math.\ Zeitschrift {\bf 195} (1987)
  197--204&

\bibitem [Pau97]&Paun, M.:& Sur le groupe fondamental des vari\'et\'es
  k\"ahl\'eriennes compactes \`a classe de Ricci num\'eriquement effective;&
  C.\ R.\ Acad.\ Sci.\ {\bf 324} (1997) 1249--1254&

\bibitem [Pau98a]&Paun, M.:& Sur les vari\'et\'es k\"ahl\'eriennes compactes
  \`a classe de Ricci num\'eriquement effective;& Bull.\ Sci.\ Math.\ {\bf 122}
  (1998) 83--92&

\bibitem [Pau98b]&Paun, M.:& Sur l'effectivit\'e num\'erique des
  images inverses de fibr\'es en droites;& Math.\ Ann.\ {\bf 310}
  (1998) 411--421&

\bibitem [Pau98c]&Paun, M.:& On the Albanese map of compact K\"ahler manifolds
  with almost positive Ricci curvature;& Pr\'epublication Institut Fourier
  Grenoble, 1998&

\bibitem [Pet93]&Peternell, Th.:& Minimal varieties with trivial
  canonical classes, I;& Math.\ Zeit\-schrift {\bf 217} (1994) 377--407&

\bibitem [Pet98]&Peternell, Th.:& Towards a Mori theory on compact K\"ahler
  $3$-folds, II;& Math.\ Ann.\ {\bf 311} (1998) 729--764&

\bibitem [Pet99]&Peternell, Th.:& Towards a Mori theory on compact K\"ahler
  $3$-folds, III;& Preprint Bayreuth Universit\"at 1999&

\bibitem [PS97]&Peternell, Th.; Serrano, F.:& Threefolds with
  numerical effective anticanonical classes;&
  Collectanea Math.\ {\bf 49} (1998) 465--517&

\bibitem [Po69]&Potters, J.:& On almost homogeneous compact complex
  surfaces;& Invent.\ Math.\ {\bf 8} (1969) 244--266&

\bibitem [Ri68]&Richberg, R.:& Stetige streng pseudokonvexe Funktionen;&
  Math.\ Ann.\ {\bf 175} (1968) 257--286&

\bibitem [Siu74]&Siu, Y.T.:& Analyticity of sets associated to Lelong
  numbers and the extension of closed positive currents;&
  Invent.\ Math.\ {\bf 27} (1974) 53--156&

\bibitem [Sk72]&Skoda, H.:& Sous-ensembles analytiques d'ordre fini ou infini
  dans $\bC^n$;& Bull.\ Soc.\ Math.\ France {\bf 100} (1972) 353--408&

\bibitem [Ta97]&Takegoshi, K.:& On cohomology groups of nef line
  bundles tensorized with multiplier ideal sheaves on compact
  K\"ahler manifolds;& Osaka J.\ Math.\ {\bf 34} (1997) 783--802&

\bibitem [Ue75]&Ueno, K.:& Classification theory of algebraic
  varieties and compact complex spaces;& Lecture Notes in Math.\
  {\bf 439}, Springer-Verlag, 1975&

\bibitem [Ue87]&Ueno, K.:& On compact analytic threefolds with
  non-trivial Albanese tori;& Math.\ Ann.\ {\bf 278} (1987) 41--70&

\bibitem [Zh96]&Zhang, Q.:& On projective manifolds with nef
  anticanonical bundle;& J.f.d.\ reine u.\ angew.\ Math.\ {\bf 478}
  (1996) 57--60&

}
\vskip30pt
\noindent
(version of October 10, 2000$\,$; printed on \today, \timeofday)

\end